\documentclass[oneside,english]{amsart}
\usepackage{mathptmx}
\usepackage{helvet}
\usepackage{courier}
\usepackage[T1]{fontenc}
\usepackage[latin9]{inputenc}
\usepackage{graphicx}
\usepackage{amssymb}

\makeatletter
\numberwithin{equation}{section} 
\numberwithin{figure}{section} 
  \theoremstyle{plain}
  \newtheorem*{thm*}{Theorem}
  \theoremstyle{plain}
  \newtheorem{thm}{Theorem}[section]
  \theoremstyle{plain}
  \newtheorem{lem}[thm]{Lemma}
  \theoremstyle{plain}
  \newtheorem{cor}[thm]{Corollary}

\usepackage{babel}
\makeatother

\begin{document}

\title{On the Minkowski Measure}

\author{Linas Vepstas <linasvepstas@gmail.com>}

\date{4 September 2008}

\begin{abstract}
The Minkowski Question Mark function relates the continued-fraction
representation of the real numbers to their binary expansion. This
function is peculiar in many ways; one is that its derivative is 'singular'.
One can show by classical techniques that its derivative must vanish
on all rationals. Since the Question Mark itself is continuous, one
concludes that the derivative must be non-zero on the irrationals,
and is thus a discontinuous-everywhere function. This derivative is
the subject of this essay.

Various results are presented here: First, a simple but formal measure-theoretic
construction of the derivative is given, making it clear that it has
a very concrete existence as a Lebesgue-Stieltjes measure, and thus
is safe to manipulate in various familiar ways. Next, an exact result
is given, expressing the measure as an infinite product of piece-wise
continuous functions, with each piece being a Möbius transform of
the form $(ax+b)/(cx+d)$. This construction is then shown to be the
Haar measure of a certain transfer operator. A general proof is given
that any transfer operator can be understood to be nothing more nor
less than a push-forward on a Banach space; such push-forwards induce
an invariant measure, the Haar measure, of which the Minkowski measure
can serve as a prototypical example. Some minor notes pertaining to
it's relation to the Gauss-Kuzmin-Wirsing operator are made.
\end{abstract}

\keywords{Minkowski Question Mark, Transfer Operator, Pushforward, Invariant
Measure.}

\subjclass[2000]{37C40 (Primary), 37C30, 37A05, 37E05, 33E99 (Secondary)}

\maketitle

\section{Introduction}

\begin{figure}
\label{cap:The-Minkowski-Question}\caption{The Minkowski Question Mark Function}
\includegraphics[bb=0bp 50bp 410bp 302bp,width=1\columnwidth]{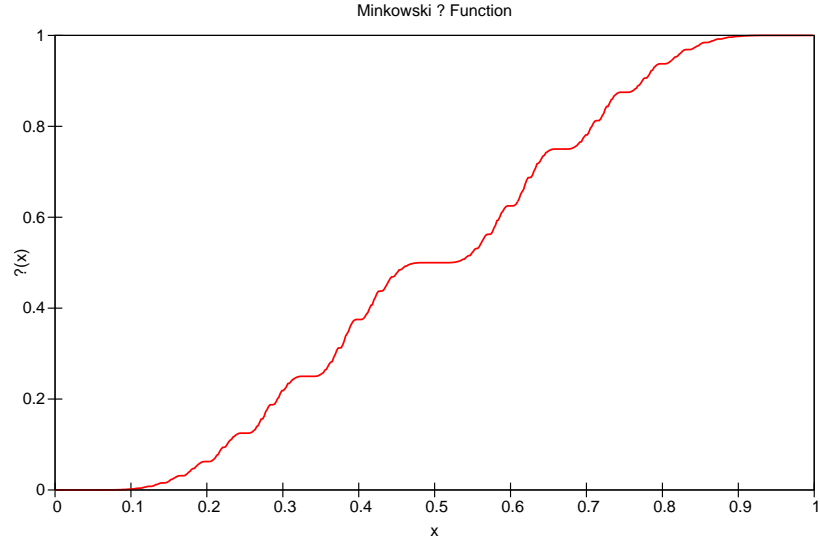}

\end{figure}

The Minkowski Question Mark function was introduced by Hermann Minkowski
in 1904 as a continuous mapping from the quadratic irrationals to
the rational numbers\cite{Mink04}. This mapping was based on the
consideration of Lagrange's criterion, that a continued fraction is
periodic if and only if it corresponds to a solution of a quadratic
equation with rational coefficients. Minkowski's presentation is brief,
but is based on a deep exploration of continued fractions. Minkowski
does provide a graph of the function; a modern version of this graph
is shown in figure \ref{cap:The-Minkowski-Question}. 

Arnaud Denjoy gave a detailed study of the function in 1938\cite{Den38},
giving both a simple recursive definition, as well as expressing it
as a summation. Its connection to the continued fractions allowed
him to provide explicit expressions for its self-similarity, described
as the action of Möbius transformations. Although Denjoy explores
a generalization that is defined for the non-negative real number
line, this essay will stick to Minkowski's original definition, defined
on the unit interval. Little is lost by sticking to the original form,
while keeping the exposition simpler.

For rational numbers, Denjoy's summation is written as \begin{equation}
?\left(x\right)=2\sum_{k=1}^{N}\left(-1\right)^{k+1}\;2^{-(a_{1}+a_{2}+...+a_{k})}\label{eq:direct ? fun defn}\end{equation}
where the $a_{k}$ are the integers appearing in the continued fraction\cite{Khin35}
expansion of $x$: \[
x=[a_{0};a_{1},a_{2},\cdots,a_{N}]=a_{0}+\frac{1}{a_{1}+\frac{1}{a_{2}+\frac{1}{a_{3}+\cdots}}}\]
For the case of general, irrational $x$, one simply takes the limit
$N\to\infty$; it is not hard to show that the resulting function
is well-defined and continuous in the classical sense of delta-epsilon
limits.

\begin{figure}
\label{fig:trees}\caption{Dyadic and Stern-Brocot Trees}
\includegraphics[bb=0bp 50bp 450bp 334bp,width=0.49\textwidth]{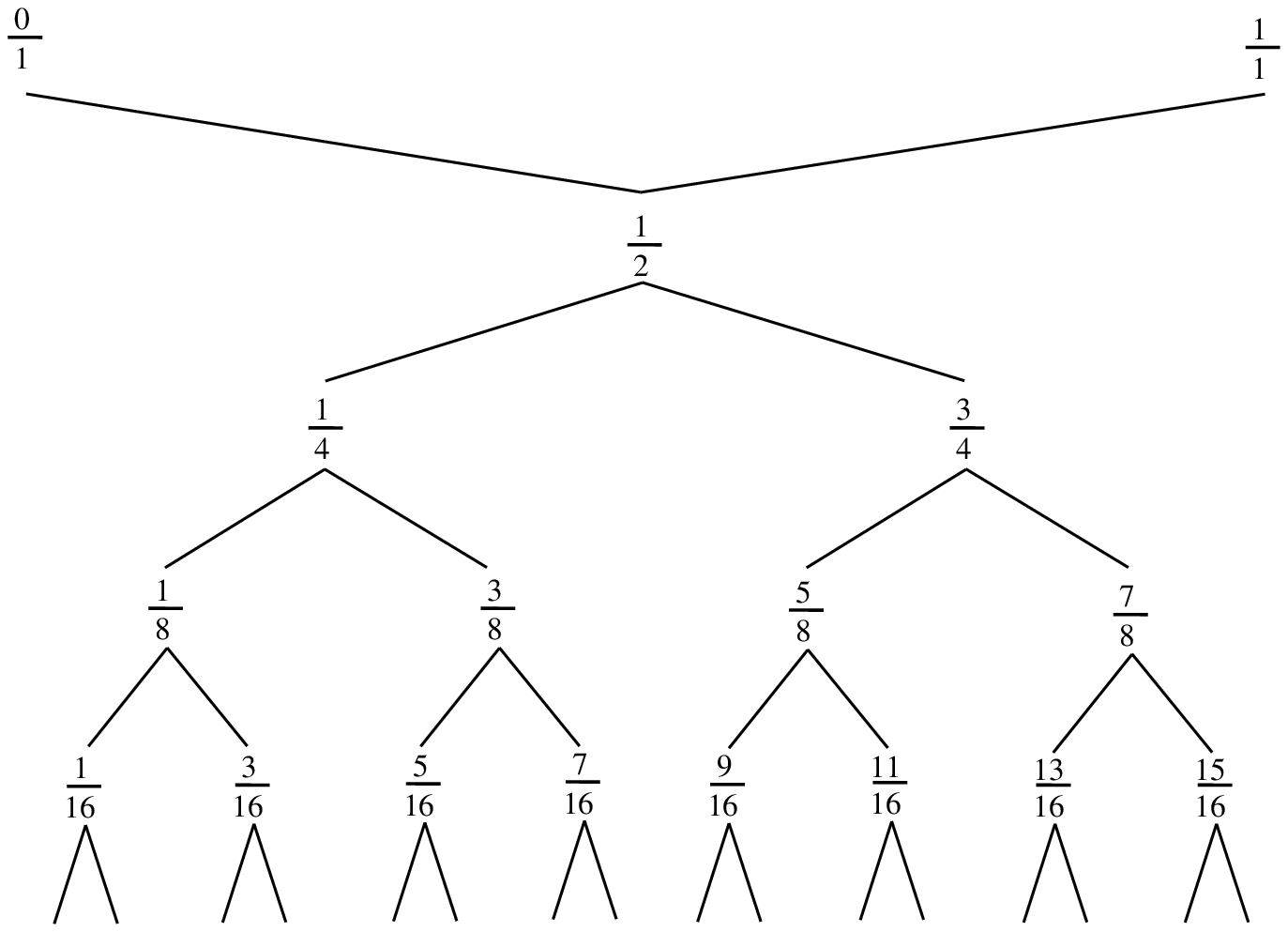}\includegraphics[bb=0bp 50bp 450bp 384bp,width=0.49\textwidth]{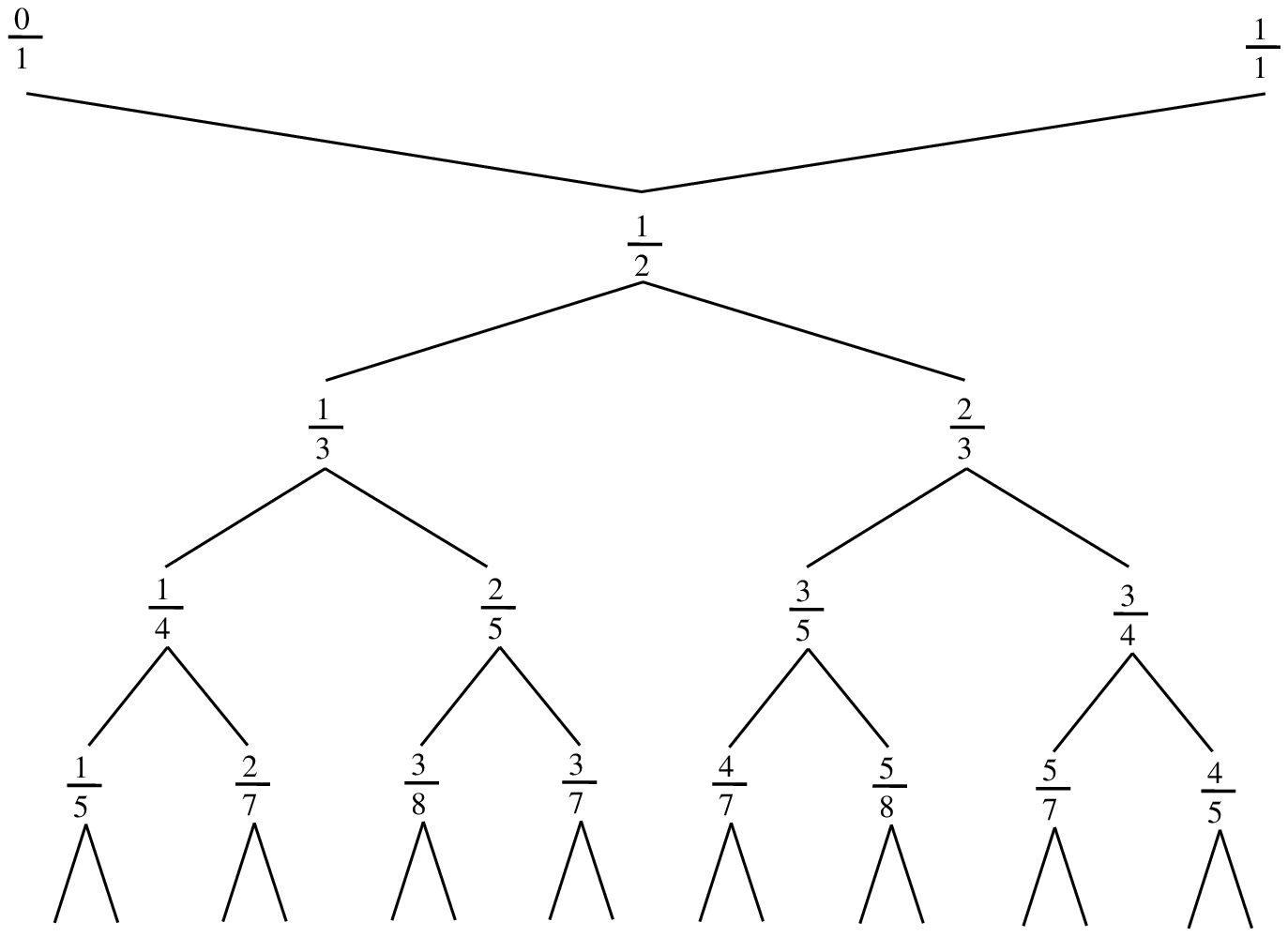}
\end{figure}

The function is symmetric in that $?\left(1-x\right)=1-?\left(x\right)$,
and it has a fractal self-similarity which is generated by \cite{Ve-M04}
\begin{equation}
\frac{1}{2}?\left(x\right)=?\left(\frac{x}{1+x}\right)\label{eq:self-similarity}\end{equation}
An equivalent, alternative definition of the function can be given
by stating that it is a correspondence between rationals appearing
in the Stern-Brocot, or Farey tree, and the binary tree of dyadic
rationals, shown in figure \ref{fig:trees}; it maps the Farey tree
into the dyadic tree.

Many insights into the nature and structure of the function, including
self-similarity and transformation properties, can be obtained by
contemplating the repercussions of this infinite binary-tree representation.
Most important is perhaps the insight into topology: the presence
of the tree indicates a product topology, and guarantees that the
Cantor set will manifest itself in various ways. 

The Question Mark has many peculiarities, and one is that it's derivative
is {}``singular'' in an unexpected way\cite{Salem1943}. One can
show, by classical techniques, that its derivative must vanish on
all rationals: it is a very very flat function, as it approaches any
rational. Since the Question Mark itself is continuous, one concludes
that the derivative must be non-zero, and infinite, on the irrationals,
and is thus a discontinuous-everywhere function. The derivative of
the Minkowski Question Mark function is interpretable as a measure,
and is often called the 'Minkowski measure'. It is this measure that
is the subject of this essay.

The derivative has been a subject of recent study. Dushistova, \emph{et
al.} show that the derivative vanishes not only on the rationals,
but also on a certain class of irrationals\cite{Dushistova2007};
this refines a result from Paridis\cite{Paridis2001}. Kesseböhmer
computes the Hausdorff dimension of the sets on which the deriviative
is non-zero\cite{Kesse2008}.

The remainder of this text is structured as follows: the next section
introduces the derivative, or 'Minkowski measure', as the statistical
distribution of the Farey numbers over the unit interval. The third
section reviews the lattice-model or cylinder-set topology for binary
strings, while the fourth presents an exact result for the Minkowski
measure in this topology. This results disproves a conjecture by Mayer\cite{May91}
that the distribution is given by the Kac model\cite{Kac1959}. In
the fifth section, the language of lattice models is converted back
to the standard topology on the real number line, to give an exact
result for the measure as an infinite product of piece-wise $C^{\infty}$
functions, with each piece taking the form of a Möbius transform $\left(ax+b\right)/\left(cx+d\right)$.
The sixth section provides a short note on modularity; the seventh
reviews an application, the eighth some generalizations. The nineth
section introduces the transfer operator, the tenth section provides
an abstract theoretical framework, demonstrating that the transfer
operator should be regarded as a pushforward. This is found to induce
an invariant measure, the Haar measure. The next section provides
a grounding for the the previous abstract discussions, illustrating
each theoretical construction with its Minkowski measure analog. The
next two sections review a twisted form of the Bernoulli operator,
and the Gauss-Kuzmin-Wirsing operator, respectively. The papaer concludes
with a call to develop machinery to better understand the discrete
eigenvalue spectrum of these operators.

\section{The Minkowski Measure }

It can be shown that the Stern-Brocot or Farey tree enumerates all
of the rationals in the unit interval exactly once\cite{Stern1858,Brocot1860}.
By considering this tree as a statistical experiment, as the source
of values for a random variable, one obtains a statistical distribution;
this distribution is exactly the Minkowski measure.

\begin{figure}
\caption{\label{cap:Distribution-of-Farey}Distribution of Farey Fractions}

\includegraphics[angle=270,width=1\columnwidth]{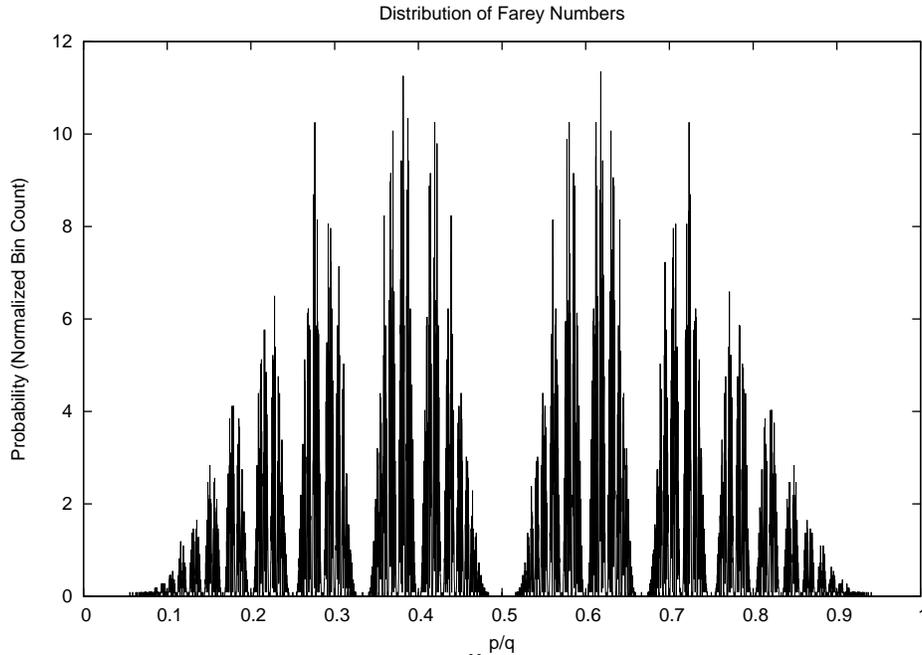}

This figure shows a histogram of the first $2^{N}=65536$ Farey fractions
according to their location on the unit interval. To create this figure,
the unit interval is divided up into $M=6000$ bins, and the bin count
$c_{m}$ for each bin $m$ is incremented whenever a fraction $p/q$
falls within the bin: $(m-1)/M\leq p/q<m/M$. At the end of binning,
the count is normalized by multiplying each bin count $c_{n}$ by
$M/2^{N}$. The normalization ensures that the histogram is of measure
one on the unit interval: that $\sum_{m=0}^{M}c_{m}=1$. 

The shape of this histogram is dependent on the choice of $M$ and
$N$. However, it does exactly capture the results of performing an
experiment on the random variable $X_{N}$ and tallying the results
into evenly-spaced bins. The result is 'exact' in the sense of eqn
\ref{eq:integral}: the height of each bin is precisely given by the
differences of the value of ? at its edges, up to the granularity
imposed by $N$. As $N\to\infty$, the errors are uniformly bounded
and vanishing. 
\end{figure}

The statistical experiment may be run as follows: consider the first
$N$ rows of the tree: these provide a pool of $2^{N}$ rationals.
A random variable $X_{N}$ is then defined by picking randomly from
this pool, with equal probability. Ultimately, the goal is to consider
the random variable $X=\lim_{N\to\infty}X_{N}$, and understand its
distribution. From the measure-theory perspective, one is interested
in the the probability $\Pr$ that a measurement of $X$ will fall
in the interval $[a,b]$ of the real number line. The distribution
is the measure $\mu\left(x\right)$ underlying this probability: \[
\Pr\left[a\le X<b\right]=\int_{a}^{b}d\mu=\int_{a}^{b}\mu\left(x\right)\, dx\]
It is straightforward to perform this experiment, for finite $N$,
on a computer, and to graph the resulting distribution. This is shown
in figure \ref{cap:Distribution-of-Farey}. What is most notable about
this figure is perhaps that it clearly belongs to a class of functions
called 'multi-fractal measures'\cite{Gut88,Man95}.

The measure $\mu$ is known as the 'Minkowski measure'; its relationship
to the Question Mark function is straight-forward, and may be stated
as a formally:

\begin{thm*}
The integral of the Stern-Brocot distribution is the Minkowski Question
Mark function; that is, \begin{equation}
?\left(b\right)-?\left(a\right)=\int_{a}^{b}d\mu\label{eq:integral}\end{equation}

\end{thm*}
\begin{proof}
The result follows from the duality of the ordinary dyadic tree and
the Stern-Brocot tree. Consider a different statistical experiment,
performed much as before, but where this time, one picks fractions
from the dyadic tree. Call the resulting random variable $Y$, understood
as the limit of discrete variables $Y_{N}$ with $N\to\infty$. Since
the dyadic rationals are uniformly distributed on the unit interval,
the corresponding measure is trivial: its just the uniform distribution.
Yet this 'trivial' experiment differs from that giving $X$ only by
the labels on the tree. Performing the same experiment, in the same
way, but just with different labellings of the values, is simply a
mapping from one set to another; the mapping is already given by eqn
\ref{eq:direct ? fun defn}, it is $?\left(x\right)$. The transformation
of measures under the action of maps is a basic result from probability
theory; the relationship between the random variables is simply that
$Y_{N}=?\left(X_{N}\right)$ and $Y=?\left(X\right)$. The measure,
of course, transforms inversely, and so by standard measure theory,
one has \[
\Pr\left[a\le X<b\right]=?(b)-?(a)\]
essentially, concluding the proof of the assertion \ref{eq:integral}.
\end{proof}
Because the equation \ref{eq:integral} resembles a classical integral,
one is sorely tempted to 'abuse the notation', and use the classical
notation to write: \begin{equation}
\frac{d?\left(x\right)}{dx}=?^{\prime}\left(x\right)=\mu\left(x\right)\label{eq:abuse}\end{equation}
The use of the derivative notation in the above is rather misleading:
properly speaking, the derivative of the Question Mark is not defined
when using the standard topology on the real number line. However,
the treatment of $\mu$ as a measure can be very rigorously founded.
Specifically, $\mu$ can be seen to be an example of an outer measure;
by the measure extension theorem (see, for example \cite{Kle08}),
an outer measure can be extended to a unique measure on the unit interval
of the real number line.

This is done as follows. Let $2=\{0,1\}$ be the set containing two
elements; let $2^{\omega}$ be the (semi-)infinite product space $2\times2\times2\times\cdots$.
This space consists of all (infinitely-long) strings in two letters.
An outer measure is a set function $\nu\,:\,2^{\omega}\to[0,\infty]$
that is monotone (\emph{i.e}. $\nu\left(A\right)<\nu\left(B\right)$
whenever $A\subset B\subset2^{\omega}$), that is sigma-subadditive
(\emph{i.e.} given any countable set $A_{1},\, A_{2},\,\ldots$ of
subsets of $2^{\omega}$ and $A\subset\bigcup_{k=0}^{\infty}A_{k}$,
then $\nu\left(A\right)\le\sum_{k=0}^{\infty}\nu\left(A_{k}\right)$
) and that $\nu$ of the empty set is zero. It is straight-forward
to show that $\nu=?^{\prime}\circ?^{-1}$ satisfies each of these
properties for an outer measure; this essentially follows from the
total ordering of the fractions on the Stern-Brocot tree (as well
as the total ordering of the dyadic rationals on the dyadic tree).
The remaining step is to notice that the {}``natural'' topology
on $2^{\omega}$ is the so-called {}``weak topology'' (sometimes
called the {}``cylinder set topology''), and that the weak topology
is finer than the natural topology on the real number line. That is,
every subset of the unit interval of the real number line is also
a subset of $2^{\omega}$. The measure extension theorem then states
that there is a single, unique measure $\tilde{\nu}$ on the unit
interval of the real number line that is identical to $\nu$ on all
subsets of the unit interval. 

It is in this way that $?^{\prime}\circ?^{-1}$ should be understood
as a measure on the unit interval, and, more precisely, a Lebesgue-Stieltjes
measure. Since $?$ is continuous and invertible (one-to-one and onto),
then $?^{\prime}$ is a Lebesgue-Stieltjes measure as well. It is
in this way that the use of the word {}``derivative'' is not entirely
inappropriate: Lebesgue-Stieltjes measures capture the essence of
integration, and insofar as the integral is an anti-derivative, one
can fairly casually talk about derivatives. In particular, later sections
below will make use of the chain rule for derivatives $\left(?\circ f\right)^{\prime}=\left(?^{\prime}\circ f\right)\cdot f^{\prime}$;
the validity of the chain rule in this context works, because it can
be seen to be a change-of-variable of a Lebesgue-Stieltjes integral. 

The construction of $?^{\prime}\circ?^{-1}$ is such that it is also
a measure (an outer measure!) on the even finer topology of $2^{\omega}$.
All these various points will be implicit in what follows, where a
more casual language will be adopted.

The product topology is also the setting for the so-called 'lattice
models' of statistical physics. The next section explicitly develops
the basic vocabulary of lattice models. The notion of fractal self-similarity
can then be explicitly equated to shifts (translations) of the lattice.
Equipped with this new vocabulary, it becomes possible to give an
exact expression for the Minkowski measure $\mu$ to be expressed
as an infinite product of Möbius functions.

\section{Lattice models}

This section provides a review of the language of lattice models,
setting up a minimal vocabulary needed for the remaining development.
The use of lattice models is partly suggested by the developments
above: the natural setting for probability theory is the product topology,
and the natural setting for the product topology is the lattice model.
The specific applicability of lattice models to the Minkowski measure
is made clear in a later section.

Consider again the set of all possible (infinite length) strings in
two letters. Each such string can be represented as $\sigma=(\sigma_{0},\sigma_{1},\sigma_{2},\ldots)$
where $\sigma_{k}\in\{-1,+1\}$. These strings can be trivially equated
to binary numbers, by writing the $k$'th binary digit as $b_{k}=(\sigma_{k}+1)/2$,
so that for each configuration $\sigma,$ one has the real number
$0\le x\le1$ given by \begin{equation}
x\left(\sigma\right)=\sum_{k=0}^{\infty}b_{k}2^{-\left(k+1\right)}\label{eq:binary-number}\end{equation}
In physics, such strings are the setting for one-dimensional lattice
models, the most famous of which is the Ising model. One considers
each number $k$ to represent a position in space, a lattice point.
At each lattice point, one imagines that there is an atom, interacting
with its neighbors, forming a perfect crystal. The goal of physics
is to describe the macroscopic properties of this crystal given only
microscopic details. There are two critical details that enter into
all such models: the principle of maximum entropy, and translational
invariance. These two principles, taken together, allow the Minkowski
measure to be specified in an exact way. 

The principle of maximum entropy states that all possible states of
the model should be equally likely, subject to obeying any constraints
that some or another average value or expectation value of a function
of the lattice be fixed. There is a unique probability distribution
that satisfies the principle of maximum entropy, it is known as the
\emph{Boltzmann distribution} or \emph{Gibbs measure}. It is given
by \begin{equation}
\Pr(\sigma)=\frac{1}{Z(\Omega)}\;\exp\left(-\beta H\left(\sigma\right)\right)\label{eq:Gibbs measure}\end{equation}
where $\Pr\left(\sigma\right)$ denotes the probability of seeing
the configuration $\sigma$. The normalizing factor $Z\left(\Omega\right)$
is known as the \emph{partition function}, it exists to ensure that
the sum of all probabilities totals unity. That is, one has \[
Z\left(\Omega\right)=\sum_{\sigma\in\Omega}\exp\left(-\beta H\left(\sigma\right)\right)\]
where $\Omega$ is the set of all possible configurations $\sigma$.
The function $H\left(\sigma\right)$ is a real-valued function, conventionally
called the 'energy' or 'Hamiltonian' in physics. The principle of
maximum entropy interprets the energy as something whose expectation
value is to be kept constant; that is, one adjusts the free parameter
$\beta$ in order to obtain a desired value for the expectation value
\[
\left\langle H\right\rangle =\frac{1}{Z\left(\Omega\right)}\sum_{\sigma\in\Omega}H\left(\sigma\right)\exp\left(-\beta H\left(\sigma\right)\right)\]

The second principle commonly demanded in physics is that of translational
invariance: that, because the crystalline lattice 'looks the same
everywhere', all results should be independent of any specific position
or location in the lattice. In particular, translational invariance
implies that $H\left(\sigma\right)$ should have the same form as
one slides around left or right on the lattice. In the present case,
the lattice is one-sided, and so one can slide in only one direction.
Defining the shift operator $\tau$ which acts on strings as $\tau(\sigma_{0},\sigma_{1},\sigma_{2},\ldots)=(\sigma_{1},\sigma_{2},\sigma_{3},\ldots)$,
the principle of translation invariance states that we should only
consider functions $H\left(\sigma\right)$ obeying the identity \[
H\left(\sigma\right)=H\left(\tau\sigma\right)\]
Unfortunately, for the present model, this cannot quite hold in the
naive sense: the lattice is not double sided (with strings of symbols
extending to infinity in both directions), but one-sided. The best
one can hope for is to isolate the part that acts just on the first
'atom' in the left-most lattice location; call this part $V\left(\sigma\right)$,
and subtract it from total: thus, instead, translational invariance
suggests that, for the one-sided lattice, one should write \begin{equation}
H\left(\sigma\right)=V\left(\sigma\right)+H\left(\tau\sigma\right)\label{eq:translational symmetry}\end{equation}
That is, the total 'energy' contained in the lattice is equal to the
energy contained in the first lattice position, plus the energy contained
in the remainder of the lattice. That the energy should be additive
(as opposed to being multiplicative, or combined in some other way)
is consistent with its use as a log-linear probability in the Gibbs
measure; that is, the use of the Gibbs measure, together with insistence
that the model be translationally invariant, implies that the energy
must be additive.

From this, it is is straight-forward to decompose $H\left(\sigma\right)$
into a sum: \begin{equation}
H\left(\sigma\right)=\sum_{k=0}^{\infty}V\left(\tau^{k}\sigma\right)\label{eq:hamilton-potential}\end{equation}
What remains is the 'potential' $V\left(\sigma\right)$, which embodies
the specifics of the problem at hand. The primary result of this note
follows from the application of the above mechanics to the Minkowski
measure: that is, to assume that the Minkowski measure is a Gibbs
measure, and to solve for the potential $V\left(\sigma\right)$.

Before doing so, there is one more relationship to make note of: the
shift operator $\tau,$ when looked at from the point of view as a
function on real numbers, looks like multiplication-by-two. To be
precise, one has the identity \[
2^{k}x\left(\sigma\right)-\left\lfloor 2^{k}x\left(\sigma\right)\right\rfloor =x\left(\tau^{k}\sigma\right)\]
where $\left\lfloor y\right\rfloor $ denotes the floor of $y$. The
left-hand side of the above is a central ingredient of the construction
of the Takagi or 'blancmange' curve; such a curve may be written as
\[
\sum_{k=0}^{\infty}\lambda^{k}f\left(2^{k}x-\left\lfloor 2^{k}x\right\rfloor \right)\]
for various functions $f$ of the unit interval. Such functions have
many curious fractal self-similarity properties, and are explored
in greater detail in \cite{Ve-T04,Ve-M04}. Writing $W=f\circ x$
shows that the Blancmange curve has the form \[
\sum_{k=0}^{\infty}\lambda^{k}W\left(\tau^{k}\sigma\right)\]
when considered on the lattice. The mapping given by eqn \ref{eq:binary-number}
is not just a mapping between strings and real numbers, it is a mapping
from the product topology to the natural topology of the real number
line. Thus, suddenly, many fractal and self-similarity phenomena are
seen to be manifestations of simple structures in the product topology,
which are then mapped to the natural topology.

\section{\label{sec:Self-Similarity}Self-Similarity }

Why should one even consider the application of lattice models to
the Minkowski measure? The core justification is that the shift operator
$\tau$ essentially acts as a multiplication-by-two, when considered
on binary the numbers $x\left(\sigma\right)$ of eqn \ref{eq:binary-number},
together with the divide-by-two self-similarity relation in eqn \ref{eq:self-similarity}.
Thus, one posits that there is a probability distribution \begin{equation}
\Pr\left(\sigma\right)=\left(?^{\prime}\circ?^{-1}\right)\left(x\left(\sigma\right)\right)\label{eq:Jacobian}\end{equation}
that is described by the Gibbs measure. The appearance here of the
inverse $?^{-1}\left(x\right)$ is so as to be able to work 'in the
same space' as where one started: that is, $?^{-1}\left(x\right)$
is a map that takes one from the space of dyadics to the space of
Farey fractions, where an operation is performed that coincidentally
takes one back to the space of dyadics. Seen in this way, eqn \ref{eq:Jacobian}
is essentially the reciprocol of a Jacobian.

In 1991, Dieter Mayer hypothesized\cite{May91} that the distribution
$\ref{eq:Jacobian}$ was given by the Kac model\cite{Kac1959}, and
specifically, by the Kac potential\begin{equation}
V(x)=\begin{cases}
-2x+1/2 & \quad\mbox{ for }0\le x\le1/2\\
2x-3/2 & \quad\mbox{ for }1/2\le x\le1\end{cases}\label{eq:Kac potential}\end{equation}
while taking $\beta=1$. The conjecture appears entirely reasonable,
at least at the numerical level: inserting this into the Gibbs measure
yields a graph that appears to be more-or-less equal the Question
Mark, up to the level of numerical detail that can be easily achieved.
However, Mayer does not provide a proof; rather, he intimates that
it must be so. 

In fact, numerical efforts suggest that there are many functions $V$
that 'come close', and can serve as the basis for similar conjectures.
Loosening the restriction of translational invariance, to allow Hamiltonians
of the form $H\left(\sigma\right)=\sum_{k=0}^{\infty}\lambda^{k}V\left(\tau^{k}\sigma\right)$
seems to allow an ever richer collection of functions V that seem
to yield good numerical approximations to the Question Mark. The Kac
potential is simply a special case; as originally formulated, it is
an interaction that decays exponentially with distance. Written in
terms of interactions between lattice points, this exponential decay
takes the form $H\left(\sigma\right)=\sum_{k=0}^{\infty}\lambda^{k}\sigma_{0}\sigma_{k}$
. Why there should be so many functions that seem to approximate the
Question Mark is not immediately apparent; perhaps it is a statement
that functions that seem 'close to' one another with respect to the
natural topology of the unit interval have dramatically different
form when considered in the product topology of the lattice model.

An exact solution has been promised; it is time to present it. To
obtain the solution, one 'follows one's nose'. Starting with\[
\Pr(x)=\left(?^{\prime}\circ?^{-1}\right)(x)\]
one clearly has $\Pr(?(x))=?^{\prime}(x)$. Applying the symmetry
from eqn \ref{eq:self-similarity}, one has that \[
\Pr\left(\frac{?(x)}{2}\right)=\Pr\left(?\left(\frac{x}{1+x}\right)\right)=?^{\prime}\left(\frac{x}{1+x}\right)\]
The division-by-two is to be recognized as the shift operator. That
is, using translational symmetry \ref{eq:translational symmetry}
in the Gibbs measure \ref{eq:Gibbs measure} yields \[
\Pr\left(\sigma\right)=\exp\left(-V\left(\sigma\right)\right)\Pr\left(\tau\sigma\right)\]
or equivalently, \[
\Pr\left(\frac{?(x)}{2}\right)=\exp\left(-V\left(\frac{?(x)}{2}\right)\right)\Pr\left(?(x)\right)\]
Back-substituting, one deduces\[
?^{\prime}\left(\frac{x}{1+x}\right)=?^{\prime}(x)\;\exp\left(-V\left(\frac{?(x)}{2}\right)\right)\]
On the other hand, one has\begin{equation}
\frac{?^{\prime}(x)}{2}=\frac{1}{(1+x)^{2}}?^{\prime}\left(\frac{x}{1+x}\right)\label{eq:meas-self-sim}\end{equation}
which may be obtained by naively differentiating the relation \ref{eq:self-similarity}.
The casual use of differentiation here deserves some remark; it is
justified by the same arguments that allow the {}``abuse of notation''
of eqn \ref{eq:abuse}; it holds as an identity involving the measure,
a change of variable taken under the integral sign. 

Continuing, the last result is plugged in to yield \[
\frac{(1+x)^{2}}{2}=\exp\left(-V\left(\frac{?(x)}{2}\right)\right)\]
which is then easily solved for $V$ to obtain

\begin{equation}
V(y)=\log2+2\log\left(1-?^{-1}(y)\right)\qquad\mbox{for }0\le y\le\frac{1}{2}\label{eq:exact-potential}\end{equation}
and $V(y)=V(1-y)$ for $1/2<x\le1$. When graphed, this potential
does appear to be vaguely tent-like, thus perhaps accounting for at
least some of the success of the Kac potential to model the Question
Mark.

To summarize the nature of this result: the potential \ref{eq:exact-potential},
let us call it the 'Minkowski potential', when used to formulate a
lattice energy \ref{eq:hamilton-potential}, gives the Minkowski measure
\ref{eq:Jacobian} as a Gibbs measure \ref{eq:Gibbs measure}. The
next section of this note deals with the form that this measure when
it is plugged back into this sequence of equations, culminating in
the Gibbs measure.

\section{A Product of Piece-wise Functions}

Plugging the exact solution \ref{eq:exact-potential} into the Hamiltonian
\ref{eq:hamilton-potential} and this in turn into the Gibbs measure
\ref{eq:Gibbs measure} allows the derivative to be written as a product
of piece-wise continuous functions. The result may be written as \begin{equation}
?^{\prime}\left(y\right)=\prod_{k=0}^{\infty}\frac{A^{\prime}\circ A_{k}\left(y\right)}{2}\label{eq:piece-product}\end{equation}
with the derivation given below. 

\begin{figure}
\label{fig:Sawteeth}\caption{Sawteeth}
 \includegraphics[width=1\columnwidth]{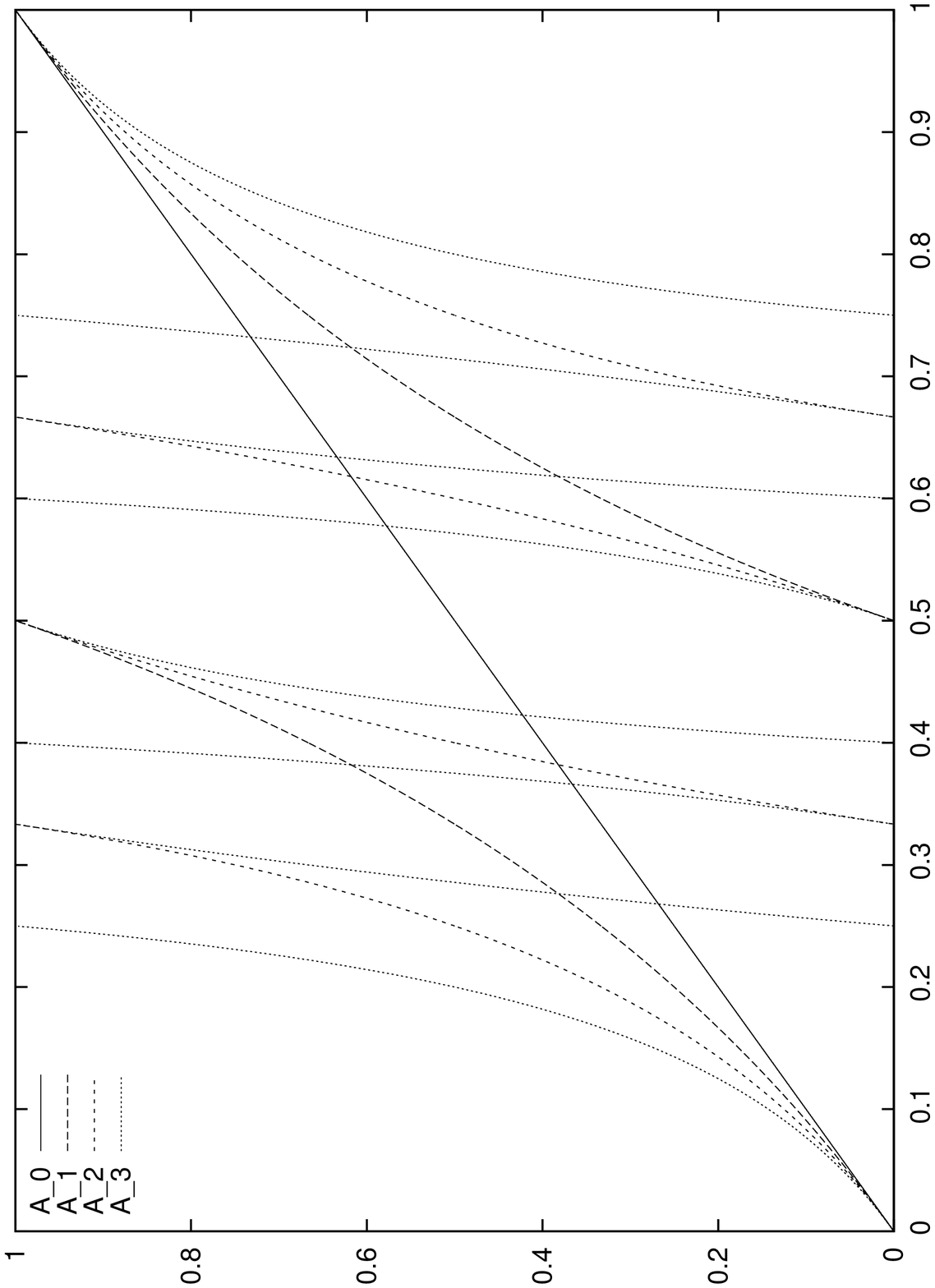}

\end{figure}

Putting together the Hamiltonian \ref{eq:hamilton-potential} with
the Gibbs measure \ref{eq:Gibbs measure} yields\begin{equation}
\Pr(\sigma)=\frac{1}{Z(\Omega)}\prod_{k=0}^{\infty}\exp\left(-V\left(\tau^{k}\sigma\right)\right)\label{eq:prod}\end{equation}
To obtain explicit expressions, it is convenient to perform the change-of-variable
$x\left(\sigma\right)=?\left(y\right)$ so that eqn \ref{eq:Jacobian}
becomes \begin{equation}
\Pr(\sigma)=?^{\prime}\left(y\right)\label{eq:measure-A}\end{equation}
while eqn \ref{eq:exact-potential} contributes 

\begin{equation}
\exp-V\left(\sigma\right)=\exp-V\left(\left(x^{-1}\circ?\right)\left(y\right)\right)=\begin{cases}
\frac{1}{2\left(1-y\right)^{2}} & \mbox{ for }0\le y\le\frac{1}{2}\\
\frac{1}{2y^{2}} & \mbox{ for }\frac{1}{2}\le y\le1\end{cases}\label{eq:exp-V}\end{equation}
The remaining terms of the product need to be evaluated by considering
the action of the shift operator $\tau$. Under the change of variable
$y=?^{-1}\left(x\left(\sigma\right)\right)$, this is given by \[
?^{-1}\left(x\left(\tau^{k}\sigma\right)\right)=?^{-1}\left(\mbox{frac}\left(2^{k}?\left(y\right)\right)\right)\]
where 'frac' is the fractional part of a number, that is, $\mbox{frac}\left(x\right)=x-\left\lfloor x\right\rfloor $.
The general expression follows by iterating on the function

\begin{equation}
A\left(y\right)=?^{-1}\left(\mbox{frac}\left(2?\left(y\right)\right)\right)=\begin{cases}
\frac{y}{1-y} & \mbox{ for }0\le y\le\frac{1}{2}\\
\frac{2y-1}{y} & \mbox{ for }\frac{1}{2}\le y\le1\end{cases}\label{eq:defn-A}\end{equation}
Define the iterated function $A_{k}$ so that \[
A_{k+1}\left(y\right)=A_{k}\circ A\left(y\right)\]
with $A_{0}\left(y\right)=y$ and $A_{1}\left(y\right)=A\left(y\right)$.
Each of the $A_{k}$ is made of $2^{k}$ Möbius transforms joined
together; so, for example, the next one is \[
A_{2}\left(y\right)=\begin{cases}
\frac{y}{1-2y} & \mbox{ for }0\le y<\frac{1}{3}\\
\frac{3y-1}{y} & \mbox{ for }\frac{1}{3}\le x<\frac{1}{2}\\
\frac{2y-1}{1-y} & \mbox{ for }\frac{1}{2}\le x<\frac{2}{3}\\
\frac{3y-2}{2y-1} & \mbox{ for }\frac{2}{3}\le x\le1\end{cases}\]
 A graph of the first four is shown in figure \ref{fig:Sawteeth}. 

The structure of the product can be elucidated by relating it to the
structure of the dyadic monoid\cite{Ve-M04,Ve-T04}. A monoid is just
a group, but without inverses; composing two elements of a monoid
simply gives another element of the monoid. The dyadic monoid is a
free monoid generated by two elements $g$ and $r$, with with $r^{2}=e$,
and no constraints on $g$. Alternately, the dyadic monoid is the
free monoid in two letters $L$ and $R$; both presentations can be
seen to correspond to a binary tree. In essence, the dyadic monoid
consists of those moves that can navigate down an infinite binary
tree. 

In this problem, the functions $g(y)=y/(1-y)$ and $r(y)=1-y$ are
the generators of the dyadic monoid. Explicitly, this is\[
A_{k+1}\left(y\right)=\begin{cases}
\left(A_{k}\circ g\right)\left(y\right)=A_{k}\left(\frac{y}{1-y}\right) & \mbox{\; for }0\le y<\frac{1}{2}\\
\left(A_{k}\circ r\circ g\circ r\right)\left(y\right)=A_{k}\left(\frac{2y-1}{y}\right) & \;\mbox{ for }\frac{1}{2}\le y\le1\end{cases}\]
Alternately, it can be convenient to work with the left $L=g$ and
right $R=r\circ g\circ r$ transformations, with locations in the
binary tree identified by navigating to the left or right, composing
$L$ and $R$. 

The functions $A_{k}(y)$ reach their extremal values 0, 1 at exactly
the Farey fractions from the $k$'th row of the Stern-Brocot tree.
These may be obtained by following the movement of the endpoints 0,1
through successive iterations of $L^{-1}$ and $R^{-1}$. These navigate
{}``up'' the tree; a point of occasional confusion can be avoided
by keeping in mind that taking the inverse reverses the order of symbols,
i.e. so that $\left(f\circ h\right)^{-1}=h^{-1}\circ f^{-1}$. For
example, the third row of the Stern-Brocot tree is given by $\frac{1}{5}=L^{-1}\circ L^{-1`}\circ L^{-1}\left(\frac{1}{2}\right)$,
$\;\frac{2}{7}=L^{-1}\circ L^{-1`}\circ R^{-1}\left(\frac{1}{2}\right)$,
$\;\frac{3}{8}=L^{-1}\circ R^{-1`}\circ L^{-1}\left(\frac{1}{2}\right)$,
$\;\frac{3}{7}=L^{-1}\circ R^{-1`}\circ R^{-1}\left(\frac{1}{2}\right)$,
$\;\frac{4}{7}=R^{-1}\circ L^{-1`}\circ L^{-1}\left(\frac{1}{2}\right)$,
and so on. Note that the order seems {}``reversed'' from the natural
order in which one might be tempted to apply the transformations $L$
and $R$. 

Of particular importance is that exponential of the potential, eqn
\ref{eq:exp-V}, is just half the derivative of $A$: \[
\exp-\left(V\circ x^{-1}\circ?\right)\left(y\right)=\frac{A^{\prime}(y)}{2}\]
This is perhaps the most remarkable result of this development. That
the Gibbs measure can be written as infinite product is no surprise.
That repeated applications of the shift operator $\tau$ can be represented
as the iteration of a function comes as no surprise either. That $A$
is the proper form for the shift operator just follows from a study
of the self-similarity properties of the Question Mark. However, the
expression for the potential $V$, as written out in eqn \ref{eq:exact-potential},
seems opaque, and without particular significance. Here, the potential
is revealed to be {}``not just some arbitrary function'', but related
to the shift itself, and, more deeply, to the inverse of the transfer
operator for the system. This is examined in greater detail below.

Several final remarks are in order. First, one has that the partition
function $Z\left(\Omega\right)=1$. This is no accident, but was guaranteed
by construction. Forcing the relationship \ref{eq:Jacobian}, together
with the fact that $?(1)=1$, ensured that all normalizing factors
that might have appeared in $Z\left(\Omega\right)$ were in fact absorbed
as an additive constant by $V$. Explicitly, this scaling leads to
the constant $\log2$ in eqn \ref{eq:exact-potential}.

Caution should be used when evaluating eqn \ref{eq:piece-product}
numerically. If preparing a graph, or numerically evaluating its integral,
the product is best terminated after $N\le\log_{2}h$ iterations,
where $h$ is the step size for the numerical integral. Using $N$
larger than this bound will result in a very noisy sampling, resulting
in numerical values that no longer resemble those of $?^{\prime}$.

\section{An Application}

Perhaps the most immediate utility of the expression \ref{eq:piece-product}
is to provide insight into the question of {}``for what values of
$x$ is $?^{\prime}\left(x\right)$ zero, and for what values is it
infinite?''; such a question is explored in \cite{Dushistova2007,Paridis2001}.
So, for example: the $A_{k}$'s can be seen to peel off binary digits,
since $A_{k}\left(y\right)=?^{-1}\left(\mbox{frac}\left(2^{k}?\left(y\right)\right)\right)$.
Thus, for any rational $p/q$, one has that $?\left(p/q\right)$ has
a finite length $N$ when expressed as a string of binary digits,
and so one immediately has that $A_{k}(p/q)=0$ when $k>N$: the infinite
product terminates, and so of course, $?(p/q)=0$.

More interesting is the case of $y$ a quadratic irrational. In this
case, $?\left(y\right)$ is a rational number, and $A_{k}\left(y\right)$,
after a finite number $M$ iterations, settles into a periodic orbit,
of length $N=2^{K}$ (of course: this is a basic result from the theory
of Pellian equations). In this case, eqn \ref{eq:piece-product} may
be written as \[
?^{\prime}\left(y\right)=\prod_{j=0}^{M-1}\frac{A^{\prime}\circ A_{k}\left(y\right)}{2}\left[\prod_{k=0}^{N-1}\frac{A^{\prime}\circ A_{k+M}\left(y\right)}{2}\right]^{\omega}\]
Clearly, the periodic part dominates. Writing \[
P\left(y\right)=\prod_{k=0}^{N-1}\frac{A^{\prime}\circ A_{k+M}\left(y\right)}{2}\]
as the contribution from the periodic part, one clearly has that $?^{\prime}\left(y\right)=0$
whenever $P\left(y\right)<1$, and $?^{\prime}\left(y\right)=\infty$
when $P\left(y\right)>1$. The accomplishment of \cite{Dushistova2007}
is to provide a precise description of these two cases, in terms of
the continued fraction expansion of $y$. Thier work would appear
to rule out any cases where $P\left(y\right)=1$.

\section{Modularity}

The self-similarity transformation of the Minkowski measure resembles
that of a modular form\cite{Apo90}. That is, given a general Möbius
transform of the form $(ax+b)/(cx+d)$, with integers $a,b,c,d$ and
$ad-bc=1$, the self-similarity of the Question Mark takes the form
\[
?\left(\frac{ax+b}{cx+d}\right)=\frac{M}{2^{N}}+\left(-1\right)^{Q}\;\frac{?\left(x\right)}{2^{K}}\]
for some integers $M,N,Q,K$, depending on the specific values of
$a,b,c,d$. Differentiating the above, one readily obtains \[
?^{\prime}\left(\frac{ax+b}{cx+d}\right)=\left(cx+d\right)^{2}\;\frac{?^{\prime}\left(x\right)}{2^{K}}\]
which is almost exactly the transformation rule for a modular form
of weight 2, and is spoiled only by the factor of $2^{K}$. This suggests
a conjecture: that there may be some holomorphic function, defined
on the upper half-plane, for which the Minkowski measure is the limit
of that function, the limit being taken a it's restriction to the
real axis. Work in this direction has been given by Alkauskas, who
notes that the generating function for the moments of the Question
Mark is holomorphic in the upper half-plane\cite{Alkauskas2007},
and gives Question Mark analogs for Maas wave forms, zeta functions
and similar number-theoretic constructs\cite{Alkauskas2007a}.

This conjecture is best illustrated numerically. Take the product
\ref{eq:piece-product}, terminating it after a finite number of terms;
then graph the reciprocal. The resulting graph fairly resembles what
one would obtain, if, for example, one graphed the Dirichlet eta function
along a horizontal line, not far above the real axis. The Dirichlet
eta is a prototypical example of a modular form; that it is also commonly
written as an infinite product is also highly suggestive. These similarities
bear further elucidation.

\section{\label{sec:Generalized-Shifts}Generalized Shifts}

Several generalizations of the Question Mark are possible. To properly
define these, it is imperative to regain a clearer view of the lattice
model methods. The ingredients to the above construction made use
of a lattice $\Omega$, whose states are encoded by $\sigma$, a lattice
shift operator $\tau$, and a function $V$ on the lattice, to define
a probability: \[
\Pr_{V}(\sigma)=\frac{1}{Z(\Omega)}\prod_{k=0}^{\infty}\exp\left(-V\left(\tau^{k}\sigma\right)\right)\]
A second important ingredient was to provide an explicit identification
of lattice arrangements $\sigma$ with real numbers $y$ in the unit
interval, given by $\sigma=x^{-1}\left(?\left(y\right)\right)$. Putting
these together lead to the form \begin{equation}
P_{f}\left(y\right)=\frac{1}{Z\left(\Omega\right)}\prod_{k=0}^{\infty}f\circ A_{k}\left(y\right)\label{eq:defn-P_f}\end{equation}
The appearance of $A_{k}$ is explicitly governed by the mapping $\sigma=x^{-1}\left(?\left(y\right)\right)$;
that is, $A_{k}$ is a representation of the shift operator $\tau^{k}$,
for this particular mapping of $\sigma$ to real values. The function
$f$ is simply $\exp-V$, that is, \[
f\left(y\right)=\exp-V\left(x^{-1}\left(?\left(y\right)\right)\right)\]
The partition function $Z\left(\Omega\right)$ is re-introduced as
a reminder that the product definition only makes sense when $f$
has been scaled appropriately; that is, $f$ is arbitrary only up
to a total scale factor; $f$ exists only as an element of a projective
space. The constraint is even stronger: inappropriately normalized
$f$ lead to the formal divergence $Z\left(\Omega\right)=\infty$;
there is only one normalization for which $Z\left(\Omega\right)$
is finite, and that furthermore, by construction, one then has $Z\left(\Omega\right)=1$.

Two directions along which to generalize are now clear: one may consider
general functions $V$, and one may consider different mappings from
the lattice space $\Omega$ to the unit interval (which is compact,
or to the positive real number line, which is not). The former is
reviewed in a later section; lets treat the later case first.

The most widely used mapping is usually the much simpler one: $x=x\left(\sigma\right)$
from eqn \ref{eq:binary-number}, instead of the {}``twisted'' mapping
$\sigma=x^{-1}\left(?\left(y\right)\right)$. The simpler mapping
is commonly used to study the Bernoulli shift\cite{Dri99,Gas92,Has92},
an entry-point into the study of exactly solvable models in symbolic
dynamics. By contrast, the {}``twisted'' mapping allowed the Question
mark to be properly formulated, thus suggesting that other maps can
lead to other interesting objects. A twisted version of the Bernoulli
Map will be presented in a later section.

\begin{figure}
\label{fig:3-adic}

\caption{The 3-adic Question Mark Function}

\includegraphics[width=1\columnwidth]{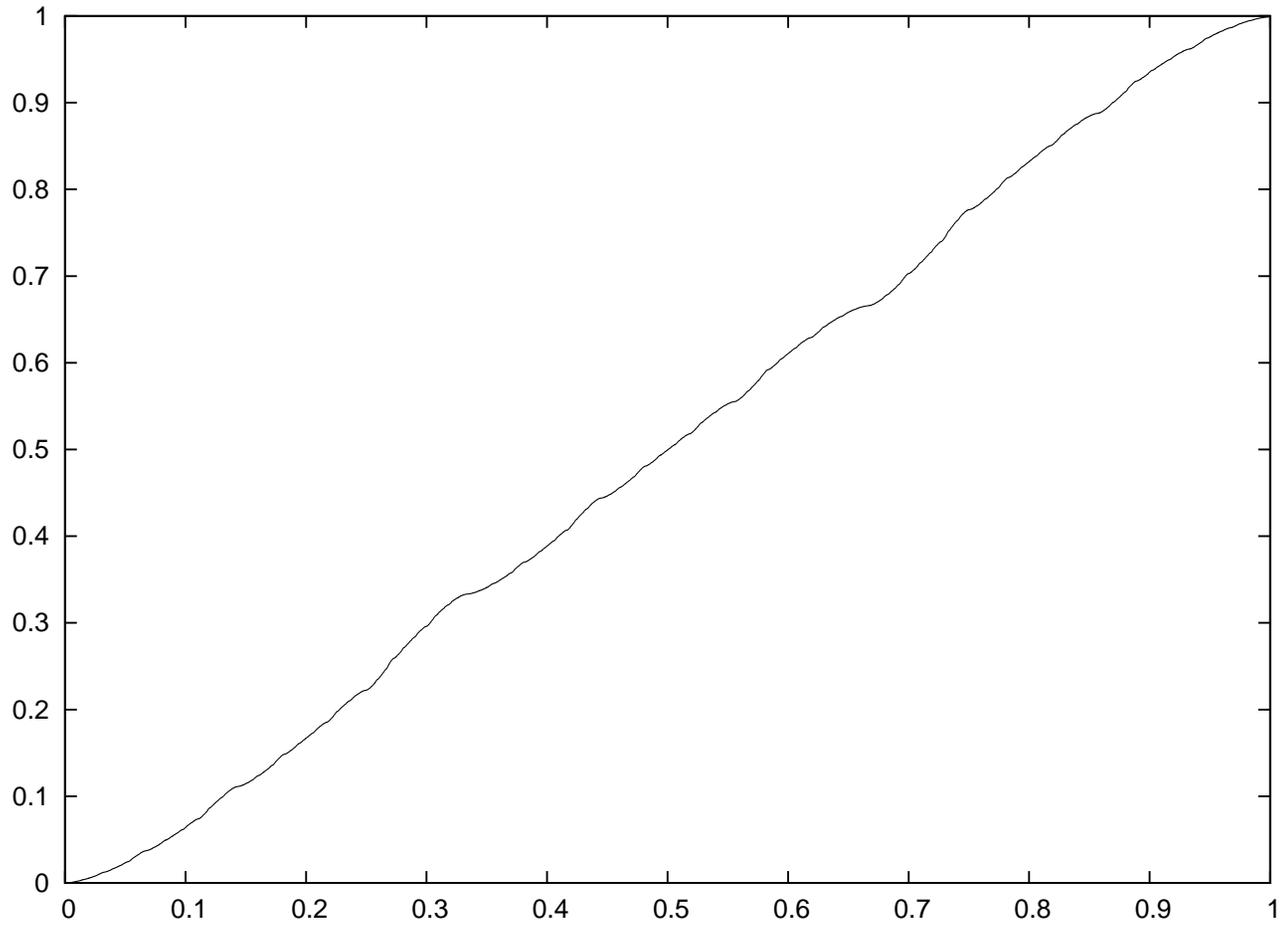}

\end{figure}

The product form \ref{eq:piece-product} provides the most concrete
avenue immediately for generalizations. One need only to consider
any differentiable, piece-wise function $A$ and use it in the construction.
For example, one may consider a 3-adic generalization of the Question
Mark, which may be constructed from \begin{equation}
B\left(y\right)=\begin{cases}
\frac{2y}{1-y} & \mbox{ for }0\le y<\frac{1}{3}\\
3y-1 & \mbox{ for }\frac{1}{3}\le y<\frac{2}{3}\\
\frac{3y-2}{y} & \mbox{ for }\frac{2}{3}\le y\le1\end{cases}\label{eq:defn-B}\end{equation}
and using $B$ instead of $A$ to generate the iterated function in
the product, and $B^{\prime}/3$ in place of $A^{\prime}/2$ as the
initial term. The resulting function, let's call it $?_{3}^{\prime}$,
can be integrated to define a 3-adic Question Mark $?_{3}$, which
is shown in figure \ref{fig:3-adic}. Note that it requires no further
normalization; one has $?_{3}(1)=1$ simply by having divided $B^{\prime}$
by the number of branches of $B$. Visually, it is considerably less
dramatic than the 2-adic Question Mark; hints of its self-similarity
are visible but less apparent.

By construction, it has self-similarity, which can be immediately
read off by considering $B^{-1}$, so, for example: \[
?_{3}\left(\frac{x}{2+x}\right)=\frac{?\left(x\right)}{3}\]
 and \[
?_{3}\left(\frac{x+1}{3}\right)=\frac{1}{3}+\frac{?\left(x\right)}{3}\]
and \[
?_{3}\left(\frac{2}{3-x}\right)=\frac{2}{3}+\frac{?\left(x\right)}{3}\]
These self-similarity properties allow an explicit equation for $?_{3}$
to be written out, entirely analogous to Denjoy's explicit summation
given in eqn \ref{eq:direct ? fun defn}, although now one must consider
even and odd values of the terms $a_{k}$ of the continued fraction
distinctly. Whether, or how, the 3-adic, or, more generally, a $p$-adic
Question Mark, inter-relate to results from Pellian equations, is
not clear. Of some curiosity is whether the feature that caused Minkowski
to present the Question Mark, namely, Lagrange's result relating periodic
continued fractions to quadratic irrationals, might somehow be reflected,
\emph{e.g}. with cubic irrationals, and an appropriate 3-adic redefinition
of periodicity in the continued fraction.

\section{The Transfer Operator}

The mapping from the lattice to the unit interval results in a collision
of effects coming from different topologies: the product topology
of the lattice, and the natural topology of the real number line.
A standard way of studying this collision is by means of the transfer
operator. The goal of this section is to introduce the transfer operator;
the next section will show, in general terms, that the transfer operator
is the push-forward of the shift operator. 

Transfer operators, also commonly called Perron-Frobenius operators
or Ruelle-Frobenius-Perron operators, as formulated by David Ruelle{[}xxx
need ref] are a commonly-used tool for studying dynamical systems\cite{Dri99}.
As the name {}``Perron-Frobenius'' suggests, these are bounded operators,
whose largest eigenvalue is equal to one; the existence of this eigenvalue
is commonly justified by appealing to the Frobenius theorem, which
states that bounded operators have a maximal eigenvalue. The origin
of the name {}``transfer'' comes from lattice-model physics, where
the {}``transfer interaction'' or the {}``transfer matrix'' is
studied\cite{Gutzwiller1990,Glimm1981,Kittel1976}. In this context,
the {}``transfer interaction'' describes the effect of neighboring
lattice sites on one-another; it can be seen to be a matrix when the
number of interacting lattice sites, and the number states a single
lattice point can be in, is finite. It becomes an operator when the
interaction is infinite range, or when the set of states at a lattice
site is no longer finite. 

It is easiest to begin with a concrete definition. Consider a function
$g:[0,1]\to[0,1]$, that is, a function mapping the unit interval
of the real number line to itself. Upon iteration, the function may
have fixed points or orbits of points. These orbits may be attractors
or repellors, or may be neutral saddle points. The action of $g$
may be ergodic or chaotic, strong-mixing or merely topologically mixing.
In any case, the language used to discuss $g$ is inherently based
on either the point-set topology of the unit interval, or the {}``natural''
topology on the unit interval, the topology of open sets. 

A shift in perspective may be gained not by considering how $g$ acts
on points or open sets, but instead by considering how $g$ acts on
distributions on the unit interval. Intuitively, one might consider
a dusting of points on the unit interval, with the local density given
by $\rho(x)$ at point $x\in[0,1]$, and then consider how this dusting
or density evolves upon iteration by $g$. This verbal description
may be given form as \begin{equation}
\rho^{\prime}(y)=\int_{0}^{1}\,\delta\left(y-g(x)\right)\rho(x)\; dx\label{eq:transfer-dirac}\end{equation}
 where $\rho^{\prime}(y)$ is the new density at point $y=g(x)$ and
$\delta$ is the Dirac delta function. 

In this viewpoint, $g$ becomes an operator that maps densities $\rho$
to other densities $\rho^{\prime}$, or notationally, \[
\mathcal{L}_{g}\rho=\rho^{\prime}\]
The operator $\mathcal{L}_{g}$ is the transfer operator or the Ruelle-Frobenius-Perron
operator. It is not hard to see that it is a linear operator, in that
\[
\mathcal{L}_{g}(a\rho_{1}+b\rho_{2})=a\mathcal{L}_{g}\rho_{1}+b\mathcal{L}_{g}\rho_{2}\]
 for constants $a,b$ and densities $\rho_{1},\rho_{2}$. 

When the function $g$ is differentiable, and doesn't have a vanishing
derivative, the integral formulation of the transfer operator \ref{eq:transfer-dirac}
can be rephrased in a more convenient form, as \begin{equation}
\left[\mathcal{L}_{g}\rho\right]\left(y\right)=\sum_{x:y=g(x)}\frac{\rho(x)}{\left|dg(x)/dx\right|}\label{eq: transfer-jacobi}\end{equation}
where the sum is presumed to extend over at most a countable number
of points. That is, when $g$ is many-to-one, the sum is over the
pre-images of the point $y$. 

Both of the above definitions are mired in a language that implicitly
assumes the real-number line, and its natural topology. Yet clearly,
this is of limited utility; extracting the topology from the definition
provides a clearer picture of the object itself.

\section{The Push-Forward}

In this section, it will be shown that the transfer operator is the
push-forward of the shift operator; a theorem and a sequence of lemmas
will be posed, that hold in general form. The point of this theorem
is to disentangle the role of topology, and specifically, the role
of measure theory, from the use of the shift operator. We begin with
a general setting.

Consider a topological space $X$, and a field $F$ over the reals
$\mathbb{R}$. Here, $F$ may be taken to be $\mathbb{R}$ itself,
or $\mathbb{C}$ or some more general field over $\mathbb{R}$. The
restriction of $F$ to being a field over the reals is is required,
so that it can be used in conjunction with a measure; measures are
always real-valued.

One may then define the algebra of functions $\mathcal{F}(X)$ on
$X$ as the set of functions $f\in\mathcal{F}(X)$ such that $f:X\to F$.
An algebra is a vector space endowed with multiplication between vectors.
The space $\mathcal{F}(X)$ is a vector space, in that given two functions
$f_{1},f_{2}\in\mathcal{F}(X)$, their linear combination $af_{1}+bf_{2}$
is also an element of $\mathcal{F}(X)$; thus $f_{1}$ and $f_{2}$
may be interpreted to be the vectors of a vector space. Multiplication
is the point-wise multiplication of function values; that is, the
product $f_{1}f_{2}$ is defined as the function $(f_{1}f_{2})(x)=f_{1}(x)\cdot f_{2}(x)$,
and so $f_{1}f_{2}$ is again an element of $\mathcal{F}(X)$. Since
one clearly has $f_{1}f_{2}=f_{2}f_{1}$, multiplication is commutative,
and so $\mathcal{F}(X)$ is also a commutative ring.

The space $\mathcal{F}(X)$ may be endowed with a topology. The coarsest
topology on $\mathcal{F}(X)$ is the \emph{weak topology}, which is
obtained by taking $\mathcal{F}(X)$ to be the space that is topological
dual to $X$. As a vector space, $\mathcal{F}(X)$ may be endowed
with a norm $\left\Vert f\right\Vert $. For example, one may take
the norm to be the $L^{p}$-norm \[
\left\Vert f\right\Vert _{p}=\left(\int\left|f(x)\right|^{p}dx\right)^{1/p}\]
For $p=2$, this norm converts the space $\mathcal{F}(X)$ into the
Hilbert space of square-integrable functions on $X$. Other norms
are possible, in which case $\mathcal{F}(X)$ has the structure of
a Banach space rather than a Hilbert space.

Consider now a homomorphism of topological spaces $g:X\to Y$. This
homomorphism induces the pullback $g^{*}:\mathcal{F}(Y)\to\mathcal{F}(X)$
on the algebra of functions, by mapping $f\mapsto g^{*}(f)=f\circ g$
so that $f\circ g:Y\to F$. The pullback is a linear operator, in
that \[
g^{*}(af_{1}+bf_{2})=ag^{*}(f_{1})+bg^{*}(f_{2})\]
That the pullback is linear is easily demonstrated by considering
how $g^{*}f$ acts at a point: $(g^{*}f)(x)=(f\circ g)(x)=f(g(x))$
and so the linearity of $g^{*}$ on $af_{1}+bf_{2}$ follows trivially.

One may construct an analogous mapping, but going in the opposite
direction, called the push-forward: $g_{*}:\mathcal{F}(X)\to\mathcal{F}(Y)$.
There are two ways of defining a push-forward. One way is to define
it in terms of the sheaves of functions on subsets of $X$ and $Y$.
The sheaf-theoretic description is more or less insensitive to the
ideas of measurability, whereas this is important to the definition
of the transfer operator, as witnessed by the appearance of the Jacobian
determinant in equation \ref{eq: transfer-jacobi}. By contrast, the
measure-theoretic push-forward captures this desirable aspect. It
may be defined as follows.

One endows the spaces $X$ and $Y$ with sigma-algebras $(X,\mathcal{A})$
and $(Y,\mathcal{B})$, so that $\mathcal{A}$ is the set of subsets
of $X$ obeying the axioms of a sigma-algebra, and similarly for $\mathcal{B}$.
A mapping $g:X\to Y$ is called {}``measurable'' if, for all Borel
sets $B\in\mathcal{B}$, one has the pre-image $g^{-1}(B)\in\mathcal{A}$
being a Borel set as well. Thus, a measurable mapping induces a push-forward
on the sigma-algebras: that is, one has a push-forward $g_{*}:\mathcal{F}(\mathcal{A})\to\mathcal{F}(\mathcal{B})$
given by $f\mapsto g_{*}(f)=f\circ g^{-1}$, which is defined by virtue
of the measurability of $g$. The push-forward is a linear operator,
in that \[
g_{*}(af_{1}+bf_{2})=ag_{*}(f_{1})+bg_{*}(f_{2})\]

One regains the transfer operator as defined in equation \ref{eq: transfer-jacobi}
by considering the limiting behavior of the push-forward on progressively
smaller sets. That is, one has

\begin{thm}
The transfer operator is the point-set topology limit of the measure-theoretic
push-forward. 
\end{thm}
\begin{proof}
The proof that follows is informal, so as to keep it simple. It is
aimed mostly at articulating the language and terminology of measure
theory. The result is none-the-less rigorous, if taken within the
confines of the definitions presented.

Introduce a measure $\mu:\mathcal{A}\to\mathbb{R}^{+}$ and analogously
$\nu:\mathcal{B}\to\mathbb{R}^{+}$. The mapping $g$ is measure-preserving
if $\nu$ is a push-forward of $\mu$, that is, if $\nu=g_{*}\mu=\mu\circ g^{-1}$.
The measure is used to rigorously define integration on $X$ and $Y$.
Elements of $\mathcal{F}(\mathcal{A})$ can be informally understood
to be integrals, in that $f(A)$ for $A\in\mathcal{A}$ may be understood
as \[
f(A)=\int_{A}\tilde{f}(z)d\mu(z)=\int_{A}\tilde{f}(z)\left|\mu^{\prime}(z)\right|dz\]
where $\left|\mu^{\prime}(x)\right|$ is to be understood as the Jacobean
determinant at a point $x\in X$. Here, $\tilde{f}$ can be understood
to be a function that is being integrated over the set $A$, whose
integral is denoted by $f(A)$. The value of $\tilde{f}$ at a point
$x\in X$ can be obtained by means of a limit. One considers a sequence
of $A\in\mathcal{A}$, each successively smaller than the last, each
containing the point $x$. One then has \[
\lim_{\overrightarrow{A\ni x}}\,\frac{f(A)}{\mu(A)}=\tilde{f}(x)\]
which can be intuitively proved by considering $A$ so small that
$\tilde{f}$ is approximately constant over $A$: \[
f(A)=\int_{A}\tilde{f}(z)d\mu(z)\approx\tilde{f}(x)\int_{A}d\mu=\tilde{f}(x)\mu(A)\]
To perform the analogous limit for the push-forward, one must consider
a point $y\in Y$ and sets $B\in\mathcal{B}$ containing $y$. In
what follows, it is now assumed that $g:X\to Y$ is a multi-sheeted
countable covering of $Y$ by $X$. By this it is meant that for any
$y$ that is not a branch-point, there is a nice neighborhood of $y$
such that its pre-image consists of the union of an at most countable
number of pair-wise disjoint sets. That is, for $y$ not a branch
point, and for $B\ni y$ sufficiently small, one may write \[
g^{-1}(B)=A_{1}\cup A_{2}\cup\cdots=\bigcup_{j=1}^{k}A_{j}\]
where $k$ is either finite or stands for $\infty$, and where $A_{i}\cap A_{j}=\varnothing$
for all $i\ne j$. At branch points, such a decomposition may not
be possible. The axiom of sigma-additivity guarantees that such multi-sheeted
covers behave just the way one expects integrals to behave: in other
words, one has \[
\mu\left(g^{-1}(B)\right)=\mu\left(\bigcup_{j=1}^{k}A_{j}\right)=\sum_{j=1}^{k}\mu\left(A_{j}\right)\]
whenever the collection of $A_{j}$ are pair-wise disjoint. Similarly,
in order to have the elements $f\in\mathcal{F}(\mathcal{A})$ behave
as one expects integrals to behave, one must restrict $\mathcal{F}(\mathcal{A})$
to contain only sigma-additive functions as well, so that \[
f\left(g^{-1}(B)\right)=f\left(\bigcup_{j=1}^{k}A_{j}\right)=\sum_{j=1}^{k}f\left(A_{j}\right)\]
As the set $B$ is taken to be smaller and smaller, the sets $A_{j}$
will become smaller as well. Denote by $x_{j}$ the corresponding
limit point of each $A_{j}$, so that $g(x_{j})=y$ and the pre-image
of $y$ consists of these points: $g^{-1}(y)=\left\{ x_{1},x_{2},\cdots\left|\, g(x_{j})=y\right.\right\} $.
One now combines these provisions to write \begin{eqnarray}
\left[g_{*}\tilde{f}\right](y) & = & \lim_{\overrightarrow{B\ni y}}\,\left[\frac{\left(g_{*}f\right)(B)}{\nu(B)}\right]\nonumber \\
 & = & \lim_{\overrightarrow{B\ni y}}\,\left[\frac{\left(f\circ g^{-1}\right)(B)}{\nu(B)}\right]\nonumber \\
 & = & \lim_{\overrightarrow{A_{j}\ni g^{-1}(y)}}\,\frac{f\left(A_{1}\cup A_{2}\cup\cdots\right)}{\nu(B)}\nonumber \\
 & = & \lim_{\overrightarrow{A_{j}\ni g^{-1}(y)}}\,\frac{\sum_{j=1}^{k}f\left(A_{j}\right)}{\nu(B)}\nonumber \\
 & = & \sum_{j=1}^{k}\tilde{f}\left(x_{j}\right)\lim_{\overrightarrow{A_{j}\ni x_{j}}}\,\frac{\mu\left(A_{j}\right)}{\nu(B)}\end{eqnarray}
The limit in the last line of this sequence of manipulations may be
interpreted in two ways, depending on whether one wants to define
the measure $\nu$ on $Y$ to be the push-forward of $\mu$, or not.
If one does take it to be the push-forward, so that $\nu=g_{*}\mu$,
then one has \[
\lim_{\overrightarrow{A_{j}\ni x_{j}}}\,\frac{\mu\left(A_{j}\right)}{g_{*}\mu(B)}=\frac{1}{\left|g^{\prime}\left(x_{j}\right)\right|}\]
where $\left|g^{\prime}\left(x_{j}\right)\right|$ is the Jacobian
determinant of $g$ at $x_{j}$. This last is a standard result of
measure theory, and can be intuitively proved by noting that $g\left(A_{j}\right)=B$,
so that \[
\nu\left(B\right)=\int_{A_{j}}g^{\prime}(z)\, d\mu(z)\approx g^{\prime}\left(x_{j}\right)\mu\left(A_{j}\right)\]
 for {}``small enough'' $B$. Assembling this with the previous
result, one has \begin{equation}
\left[g_{*}\tilde{f}\right](y)=\sum_{x_{j}\in g^{-1}(y)}\frac{\tilde{f}\left(x_{j}\right)}{\left|g^{\prime}\left(x_{j}\right)\right|}\label{eq:Transfer-point}\end{equation}
which may be easily recognized as equation \ref{eq: transfer-jacobi}.
This concludes the proof of the theorem, that the transfer operator
is just the point-set topology limit of the push-forward. 
\end{proof}
In simplistic terms, the push-forward can be thought of as a kind
of change-of-variable. Thus, one should not be surprised by the following
lemma, which should be recognizable as the Jacobian, from basic calculus. 

\begin{lem}
(Jacobian) One has \[
\sum_{x_{j}\in g^{-1}(y)}\frac{1}{\left|g^{\prime}\left(x_{j}\right)\right|}=1\]
 
\end{lem}
\begin{proof}
This follows by taking the limit $\overrightarrow{A_{j}\ni x_{j}}$
of \[
\frac{\mu\left(A_{j}\right)}{g_{*}\mu(B)}=\frac{\mu\left(A_{j}\right)}{\sum_{i=1}^{k}\mu\left(A_{i}\right)}\]
 and then summing over $j$.
\end{proof}
\begin{cor}
\label{cor:unit-eigenvec}The uniform distribution (\emph{i.e.} the
measure) is an eigenvector of the transfer operator, associated with
the eigenvalue one. 
\end{cor}
\begin{proof}
This may be proved in two ways. From the viewpoint of point-sets,
one simply takes $\tilde{f}=\mbox{const.}$ in equation \ref{eq:Transfer-point},
and applies the lemma above. From the viewpoint of the sigma-algebra,
this is nothing more than a rephrasing of the starting point, that
$\nu=g_{*}\mu$, and then taking the space $Y=X$, so that the push-forward
induced by $g:X\to X$ is a measure-preserving map: $g_{*}\mu=\mu$. 
\end{proof}
The last corollary is more enlightening when it is turned on its side;
it implies two well-known theorems, which follow easily in this framework.

\begin{cor}
(Ruelle-Perron-Frobenius theorem). All transfer operators are continuous,
compact, bounded operators; furthermore, they are isometries of Banach
spaces.
\end{cor}
\begin{proof}
This theorem is of course just the Frobenius-Perron theorem, recast
in the context of measure theory. By definition, the measures have
unit norm: that is, $\left\Vert \mu\right\Vert _{1}=1$ and $\left\Vert \nu\right\Vert _{1}=1$.
This is nothing more than the statement that the spaces $X$ and $Y$
are measurable: the total volume of $X$ and $Y$ is, by definition,
one. Since $\nu=g_{*}\mu$, we have $\left\Vert g_{*}\mu\right\Vert _{1}=1$,
and this holds for all possible measures $\mu\in\mathcal{F}(X)$. 

Recall the definition of a bounded operator. Given a linear map $T:U\to V$
between Banach spaces $U$and $V$, then $T$ is bounded if there
exists a constant $C<\infty$ such that $\left\Vert Tu\right\Vert _{V}\le C\left\Vert u\right\Vert _{U}$for
all $u\in U$. But this is exactly the case above, with $T=g_{*}$,
and $U=\mathcal{F}(X)$, $V=\mathcal{F}(Y)$, and $C=1$. The norm
of a bounded operator is conventionally defined as\[
\left\Vert T\right\Vert =\sup_{u\ne0}\frac{\left\Vert Tu\right\Vert _{V}}{\left\Vert u\right\Vert _{U}}=\sup_{\left\Vert u\right\Vert _{U}\le1}\left\Vert Tu\right\Vert _{V}\]
and so we have the norm of $g_{*}$ being $\left\Vert g_{*}\right\Vert =1$.
That $g_{*}$ is an isometry follows trivially from $\left\Vert g_{*}\mu\right\Vert _{1}=\left\Vert \nu\right\Vert _{1}$
and that $g_{*}$ is linear.
\end{proof}
The corollary \ref{cor:unit-eigenvec} can also be treated as a corollary
to the Perron-Frobenius theorem: namely, that there is at least one
vector that corresponds to the maximum eigenvalue of $g_{*}$. This
eigenvector is in fact the Haar measure, as the next theorem shows. 

\begin{thm}
\label{thm:Haar-measure}(Haar measure) For any homomorphism $g:X\to X$,
one may find a measure $\mu$ such that $g_{*}\mu=\mu$; that is,
every homomorphism $g$ of $X$ induces a measure $\mu$ on $X$ such
that $g$ is a measure-preserving map. If $g$ is ergodic, then the
measure is unique. 
\end{thm}
\begin{proof}
By definition, $\mu$ is a fixed point of $g_{*}$. The fixed point
exists because $g_{*}$is a bounded operator, and the space of measures
is compact, and so a bounded operator on a compact space will have
a fixed point. The existence of the fixed point is given by the Markov-Kakutani
theorem\cite[p 456]{Dunford1957}. The Markov-Kakutani theorem also
provides the uniqueness condition: if there are other push-forwards
$h_{*}$ that commute with $g_{*}$, then each such push-forward will
also have a fixed point. The goal is then to show that when $g$ is
ergodic, there are no other functions $h$ that commute with $g$.
But this follows from the definition of ergodicity: when $g$ is ergodic,
there are no invariant subspaces, and the orbit of $g$ is the whole
space. As there are no invariant subspaces, there are no operators
that can map between these subspaces, \emph{i.e.} there are no other
commuting operators.
\end{proof}
A peculiar special case is worth mentioning: if $g$ is not ergodic
on the whole space, then typically one has that the orbit of $g$
splits or foliates the measure space into a bunch of pairwise disjoint
leaves, with $g$ being ergodic on each leaf. The Markov-Kakutani
theorem then implies that there is a distinct fixed point $\mu$ in
each leaf, and that there is a mapping that takes $\mu$ in one leaf
to that in another. 

In the language of dynamical systems, the push-forward $g_{*}$ is
commonly written as $\mathcal{L}_{g}$, so that one has \[
g_{*}=\mathcal{L}_{g}\,:\,\mathcal{F}\left(X\right)\to\mathcal{F}\left(X\right)\]
now being called the transfer operator or the Ruelle-Frobenius-Perron
operator.

In the language of physics, the fixed point $\mu$ is called the {}``ground
state'' of a system. When it is unique, then the ground state is
not degenerate; when it is not unique, then the ground state is said
to be degenerate. The operator $g_{*}$ is the time-evolution operator
of the system; it shows how physical fields $f\in\mathcal{F}(X)$
over a space $X$ evolve over time. When $F$ is the complex numbers
$\mathbb{C}$, the fact that $\left\Vert g_{*}\right\Vert =1$ is
essentially a way of stating that time-evolution is unitary; the Frobenius-Perron
operator is the unitary time-evolution operator of the system. What
is called {}``second quantization'' in physics should be interpreted
as the fitting of the space $\mathcal{F}(X)$ with a set of basis
vectors, together with a formulation of $g_{*}$ in terms of that
basis.

\section{Grounding}

To make sense of the discussion above, it is best to ground the symbols
in concrete terms. For the most part, this paper is concerned either
with the space $X=\Omega$, the space of strings in two letters, or
the spaces $X$ that are somehow isomorphic to $\Omega$. This space
is a (semi-)infinite product of $\mathbb{Z}_{2}=\{0,1\}$; elements
$\sigma\in\Omega$ of this space were values $\sigma=(\sigma_{0},\sigma_{1},\cdots)$
with each $\sigma_{k}\in\mathbb{Z}_{2}$. The function $g$ is to
be equated with the shift operator $\tau$; that is, when $X=\Omega$,
then $g=\tau$ and $g:X\to X$ is just $\tau:\Omega\to\Omega$ acting
such that $\tau(\sigma_{0},\sigma_{1},\sigma_{2},\cdots)=(\sigma_{1},\sigma_{2},\cdots)$.
Curiously, the structure of the space $\Omega$ can be understood
to be a binary tree: starting at $\sigma_{0}$ one can branch to the
left or to the right, to get to $\sigma_{1}$, and so on. 

However, a product space such as this also allows the language of
probability theory come into play, and this is where things get interesting.
Product spaces are the natural setting for probability: the element
$\sigma\in\Omega$ is called a {}``random variable'', and successive
values $\sigma_{k}$ of $\sigma$ are {}``measurements'' of the
random variable. This is also the manner in which measure theory gets
its foot in the door.

The mappings $x:\Omega\to\mathbb{R}$ given by eqn \ref{eq:binary-number}
is now seen to be an example of an element of the algebra $\mathcal{F}\left(\Omega\right)$
introduced in the previous section. The mapping $y:\Omega\to\mathbb{R}$
given by $\sigma\mapsto y=?^{-1}\left(x\left(\sigma\right)\right)$
is another example of an element of $\mathcal{F}\left(\Omega\right)$.
Both the maps $x$ and $y$ are invertible on the closed unit interval,
establishing a $1-1$ mapping between $\Omega$ and the closed unit
interval. Consider, then, the function \[
b=x\circ\tau\circ x^{-1}\]
which maps the unit interval into itself. Let $w\in[0,1]$ be a point
in the unit interval; then it is not hard to see that $b(w)=2w\mbox{ mod}1$.
In the literature, this function is commonly called the Bernoulli
shift; its corresponding random process is known as the Bernoulli
process. It is well-known to be ergodic. The other function that was
considered here was \[
A=y\circ\tau\circ y^{-1}\]
Its action on the unit interval has already been given by eqn \ref{eq:defn-A}.
In a previous section, we have already considered generalizations
of $A$, that is, other functions $B:[0,1]\to[0,1]$, such as, for
example, eqn \ref{eq:defn-B}. 

Each of these maps $b$, $A$, $B$ can be taken as an example of
a map $g:X\to X$, now taking $X$ to the unit interval. By the previous
theorem, each such map $g$ induces a push-forward $g_{*}=\mathcal{L}_{g}$,
and this push-forward then induces an invariant measure $\mu_{g}$
on $X$. For the Bernoulli map $b$, the push-forward $\mathcal{L}_{b}$
is introduced in \cite{Gas92,Has92}, and a detailed study of its
eigenvectors is given in \cite{Dri99}. The invariant measure $\mu_{b}$
is well-known to be the identity on the unit interval. This just means
that $b$ considered as an ergodic process, is uniformly distributed
on the unit interval. 

The invariant measure $\mu_{A}$ has already been established in eqn
\ref{eq:measure-A} as being the Minkowski measure; that is, $\mu_{A}=?^{\prime}$.
A previous section already gave specific mechanics for how to take
a general, differentiable onto-map $B$ and define its associated
measure $\mu_{B}$ as the generalization of eqn \ref{eq:piece-product}.
The next section introduces the transfer operator $\mathcal{L}_{A}$
corresponding to $A$.

\section{Flightiness}

To counterbalance the grounding, a brief bit of flightiness as to
this process is offered up. The map $g:X\to X$, can be understood
to give the (discrete) time-evolution of a dynamical system. It is
curious to ponder that iterating on $g$ can be thought of as the
(semi-)infinite chain \[
X\begin{array}[b]{c}
g\\
\longrightarrow\end{array}X\begin{array}[b]{c}
g\\
\longrightarrow\end{array}X\begin{array}[b]{c}
g\\
\longrightarrow\end{array}X\begin{array}[b]{c}
g\\
\longrightarrow\end{array}\cdots\]
and so any particular $g$, together with an initial state $x_{0}$,
can be interpreted as selecting a particular sequence $(x_{0},x_{1},x_{2},\ldots)$
out of the (semi-)infinite product space $\Omega=X\times X\times X\times\cdots$.
The set membership function now defines a measure on $\Omega$: sets
that contain sequences $(x_{0},x_{1},x_{2},\ldots)$ for which $x_{k+1}=g\left(x_{k}\right)$
will have a non-zero measure, all others will have a measure of zero. 

This now begs the question: what sort of operators might have this
measure as their invariant measure? Which of these operators might
be push-forwards? What would they be a push-forward of? Last but not
least, all of the spaces under discussion are metric spaces, and so
it would be curious to know {}``what else is close by''. These questions
are not entirely idle. The product space $\Omega$ can be seen to
be {}``space-time''; by analogy to the binary tree, the structure
of this space-time is once-again seen to be branching; infinitely
branching in this case. But such branching is already seen in physics:
this is nothing other than the many-worlds interpretation of quantum
mechanics, where, at each moment of time, the universe is seen to
split into an {}``infinite number of copies of itself''. The many-worlds
branching structure is not well-explored or well-understood, although
it can illuminate problems of quantum measurement and wave-function
collapse. One such curious problem is the so-called {}``Mott problem''.
The question that is posed by the Mott problem is {}``why does the
spherically symmetric wave function of an emitted alpha-ray always
leave a straight-line track in a Wilson cloud chamber?'' The answer
was given by Heisenberg and Mott in 1929; it can be re-interpreted
in the context of the current formulation as a statement that straight-line
tracks have a significant measure in this branching space-time, as
they all correspond to sets that contain the point that is the spherical
wave-function, and no other sets do. But this is just an impressionistic
daydream, lacking rigor; the hope is that it motivates further study.

\section{The Twisted Bernoulli Operator}

We know consider the transfer operator or push-forward of the mapping
$A$. It resembles the Bernoulli shift, but is twisted in a curious
way; it is a continued-fraction analog of the Bernoulli shift. It
is straight-forward to write it down; it is very simply \[
\left[\mathcal{L}_{A}f\right]\left(y\right)=\frac{1}{\left(1+y\right)^{2}}f\left(\frac{y}{1+y}\right)+\frac{1}{\left(2-y\right)^{2}}f\left(\frac{1}{2-y}\right)\]
This expression is easily obtained by starting with eqn \ref{eq: transfer-jacobi}
and substituting in $A$ from eqn \ref{eq:defn-A}.

Of course, the Minkowski measure is an eigenfunction of the twisted
Bernoulli shift; that is, one has \[
\mathcal{L}_{A}?^{\prime}=?^{\prime}\]
This can be seen in several ways: The Minkowski measure was constructed
as the Haar measure induced by $A$, and so by theorem \ref{thm:Haar-measure},
of course this identity must hold. But it can also be validated by
more pedestrian methods: simply plug in the measure, and turn the
crank, making use of the self-similarity relation \ref{eq:meas-self-sim}
as needed. 

Other eigenvalues of $\mathcal{L}_{A}$ are readily found. Consider
the functions $P_{f}$ as previously defined in \ref{eq:defn-P_f}.
The action of the twisted Bernoulli shift on $P_{f}$ is readily computed,
and is given by \[
\mathcal{L}_{A}P_{f}=P_{f}\cdot\mathcal{L}_{A}f\]
We've already observed that when $f=A^{\prime}/2$, that one has \[
\mathcal{L}_{A}\frac{A^{\prime}}{2}=1\]
so that $P_{A^{\prime}/2}=?^{\prime}$. Another, linearly independent
function which satisfies $\mathcal{L}_{A}f=1$ is \[
C\left(y\right)=\begin{cases}
\frac{1}{3\left(1-y\right)^{3}} & \mbox{ for }0\le y\le\frac{1}{2}\\
\frac{1}{3y^{3}} & \mbox{ for }\frac{1}{2}\le y\le1\end{cases}\]
 that is, $\mathcal{L}_{A}C=1$. The collection of all $f$ which
satisfy $\mathcal{L}_{A}f=1$ is large, but is easily specified. Let\[
f\left(y\right)=\begin{cases}
f_{0}\left(y\right) & \mbox{ for }0\le y\le\frac{1}{2}\\
f_{1}\left(y\right) & \mbox{ for }\frac{1}{2}\le y\le1\end{cases}\]
then, given any, arbitrary $f_{0}$, solving for $\mathcal{L}_{A}f=1$
results in \[
f_{1}\left(y\right)=\frac{1}{y^{2}}-\frac{1}{\left(3y-1\right)^{2}}f_{0}\left(\frac{2x-1}{3x-1}\right)\]
Such $f$ are properly scaled; that is, $Z\left(\Omega\right)=1$.
Because of the projective relationship between the functions $f$
and the normalization $Z\left(\Omega\right)$, it appears that this
technique will result in only in eigenfunctions $P_{f}$ with eigenvalue
1. 

Keep in mind that while $P_{f+h}\ne P_{f}+P_{h}$, such a relationship
does hold in the logarithmic sense: $P_{fh}=P_{f}P_{h}$, thus $P$
can be construed to be a (logarithmic) homomorphism of Banach spaces.

\section{The Gauss-Kuzmin-Wirsing Operator}

The Gauss-Kuzmin-Wirsing (GKW) operator is the transfer operator of
the Gauss map $h(x)=1/x-\left\lfloor 1/x\right\rfloor $. Denoting
this operator as $\mathcal{L}_{G}$, it acts on a function $f$ as
\begin{equation}
\left[\mathcal{L}_{G}f\right](x)=\sum_{n=1}^{\infty}\frac{1}{\left(x+n\right)^{2}}\, f\left(\frac{1}{x+n}\right)\label{eq:gkw}\end{equation}
It has a well-known eigenfunction, noted by Gauss, corresponding to
eigenvalue one, given by $1/(1+y)$. It is also relatively straightforward
to observe that the Minkowski measure is also an eigenfunction, also
corresponding to eigenvalue one. This is readily seen by differentiating
the Question Mark: one has \[
?\left(\frac{1}{n+x}\right)=\frac{1}{2^{n-1}}-\frac{?\left(x\right)}{2^{n}}\]
and so also \[
?^{\prime}\left(\frac{1}{n+x}\right)=\frac{\left(x+n\right)^{2}}{2^{n}}?^{\prime}\left(x\right)\]
This can be directly inserted into eqn \ref{eq:gkw} to very easily
find that \[
\mathcal{L}_{G}?^{\prime}=?^{\prime}\]
This result immediately suggests that the lattice model methods can
be applied to find other eigenfunctions. (The 3-adic Question mark,
at least as defined above, does not have this property, suggesting,
perhaps, that a different 3-adic Question Mark should be formulated).

The GKW operator is tied by various identities to the twisted Bernoulli
operator $\mathcal{L}_{A}$. These are explicitly given by \cite{Iso03};
they can be crudely summarized by saying that these identities are
the operator analogues of Denjoy's explicit summation \ref{eq:direct ? fun defn}.
In these rough terms, the Bernoulli shift chops off one digit of a
binary number at a time; these may be in turn re-assembled to form
the GKW operator.

The action of the GKW operator $\mathcal{L}_{G}$ on $P_{f}$ is not
at all tidy; the result that $\mathcal{L}_{G}?^{\prime}=?^{\prime}$
appears to be entirely due to the special form of $f=A^{\prime}/2$.
This is perhaps no surprise.

Other fractal eigenvectors for the GKW can be constructed by a technique
entirely parallel to the construction of the $P_{f}$; one only has
to use a different lattice model. For this model, construct \[
R_{f}\left(y\right)=\frac{1}{Z\left(\Omega\right)}\prod_{k=0}^{\infty}f\circ h_{k}\left(y\right)\]
where $h_{k+1}=h_{k}\circ h$ is function iteration, and $h(y)=\mbox{frac}\left(1/y\right)$
is the Gauss map, as before. This multi-branched $h$ has the property
that \[
h\left(\frac{1}{x+n}\right)=x\]
for all integer values of $n\ge1$. Applying the GKW operator to $R_{f}$
results in \[
\mathcal{L}_{G}R_{f}=R_{f}\mathcal{L}_{G}f\]
in analogy to the twisted Bernoulli operator. As before, functions
$f$ for which $\mathcal{L}_{G}f=1$ also implies that $Z\left(\Omega\right)=1$.
The means for finding solutions to $\mathcal{L}_{G}f=1$ proceeds
analogously to the twisted-Bernoulli case. Thus, let \[
f\left(y\right)=\begin{cases}
f_{1}\left(y\right) & \mbox{ for }\frac{1}{2}<y\le1\\
\cdots & \cdots\\
f_{n}\left(y\right) & \mbox{ for }\frac{1}{n+1}<y\le\frac{1}{n}\\
\cdots & \cdots\end{cases}\]
The equation $\mathcal{L}_{G}f=1$ establishes only a single constraint
between these $f_{n}$'s: all but one may be freely chosen. As an
example, consider \[
f_{n}\left(y\right)=\frac{a_{n}}{y^{2}}\]
Then the resulting $f$ satisfies $\mathcal{L}_{G}f=1$ provided that
the constants $a_{n}$ satisfy \[
\sum_{n=1}^{\infty}a_{n}=1\]
By taking $a_{n}=2^{-n}$, one regains the Question Mark function;
that is, $R_{f}=?^{\prime}$ for this case. Quick numerical work shows
that other geometric series return commonly-seen {}``Question-Mark-like''
distributions, while arithmetic series, such as $a_{n}=6/\left(\pi n\right)^{2}$,
yeild more unusual shapes.

This construction shows that the space of fractal eigenfunctions of
the GKW operator is very large; but, as before, all of these eigenfunctions
appear to be non-decaying, corresponding to eigenvalue 1.

\section{Conclusion}

The principal results of this note are the statistical distribution
of equation \ref{eq:integral}, which allows numerical explorations
of the binary trees to be explicitly equated to the Minkowski measure,
and equation \ref{eq:piece-product}, which shows that the seemingly
intractable Minkowski measure can be built up out of analytic functions
in a regular way. A sufficient amount of formal theory is developed
to show how to root this measure in a foundational apparatus, thus
enabling various basic manipulations on it to be safely performed.
Of these, the most notable was to expose the transfer operator as
a push-forward on a Banach space, leading to an interpretation of
the Minkowski measure as a invariant measure, a Haar measure, induced
by a dynamical map.

An important theme is the tension between the natural topology on
the real number line, and the much finer 'weak topology' on the product
space $2^{\omega}$. It should be clear that the fine topology is
the source for a broad variety for fractal, self-similar phenomena.
The use of the measure extension theorem allows these fine-grained
fractal phenomena to manifest themselves onto the natural topology;
the Minkowski measure serves as a non-trivial yet still simple example
of such a manifestation.

However, various open problems remain. There have now been many generalizations
of the Question Mark proposed; besides that of Denjoy, and those of
this note, some additional ones are given in \cite{Alkauskas2008}.
How are all of these related? How are they best enumerated? Is there
some characteristic labelling that best distinguishes each form? Insofar
as these variants do occasionally manifest in practical problems,
it is useful to have a coherent naming system to describe them.

It is noted in section \ref{sec:Self-Similarity} that there seem
to be many different functions $V$ that generate functions $Q$ that
seem close to the Question Mark; in some cases, these appear to be
arbitrarily close, when studied numerically. Why is this? Is there
an actual degeneracy, with different $V$'s generating the same function
$Q$? If not, then the metric on the Banach space of different $V$'s
must be wildly folded to yeild close-by $Q$'s: the metrics are not
{}``compatible'', or presumably discontinuous in some way. The structure
of the mapping from the space of $V$'s to the space of $Q$'s is
opaque, and would benefit from clarification. 

Section \ref{sec:Generalized-Shifts} considered a 3-adic form of
the Question Mark. Its construction is somewhat \emph{ad hoc}, and
it's properties are less than pretty: for example, it is not an eigenfunction
of the GKW operator. Is there another construction that would be an
eigenfunction? If not, why not? If so, what would the correct $p$-adic
generalization be? Let's assume that unambiguous $p$-adic generalizations
can be written. What are thier connections to number theory? That
is, $?\left(x\right)$ maps rationals to quadratic surds; is there
an analogous statement for such $?_{p}\left(x\right)$?

The most important open problem surrounds the observation that transfer
operators have a discrete spectrum, when the space of functions is
restricted to a much smaller set, for instance, when it is restricted
to the set of square-integrable functions, or to the set of polynomials.
In some cases, this spectrum can be obtained exactly; this has been
done for the case of the Bernoulli operator \cite{Gas92,Has92,Dri99}.
More generally, obtaining this spectrum seems almost intractable;
here, the prototypical example is the GKW operator. Of course, both
the eigenvalues and the eigenvectors of the GKW operator can be computed
numerically; however, there is no articulation or framework within
which to talk about them and to manipulate them; they remain opaque.
They seem to be mere holomorphic functions (by definition!) and oscillatory,
of course; but there is not yet a way of writing down their power
series directly.

Thus, the open problem is to develop a technical framework that will
identify the discrete eigenvalue spectrum of a transfer operator,
and to directly provide expressions for the holomorphic functions
which are their eigenvectors. Such a framework would be of broad utility. 

\bibliographystyle{plain}
\bibliography{/home/linas/src/fractal/paper/fractal}

\begin{thebibliography}{10}

\bibitem{Alkauskas2007a}
Giedrius Alkauskas.
\newblock Generating and zeta functions, structure, spectral and analytic
  properties of the moments of minkowski question mark function authors:.
\newblock {\em ArXiv}, 0801.0056, 2007.

\bibitem{Alkauskas2007}
Giedrius Alkauskas.
\newblock The moments of minkowski ?(x) function: dyadic period functions.
\newblock {\em ArXiv}, 0801.0051, 2007.

\bibitem{Alkauskas2008}
Giedrius Alkauskas.
\newblock Minkowski question mark function and its generalizations, associated
  with p-continued fractions: fractals, explicit series for the dyadic period
  function and moments.
\newblock {\em ArXiv}, 0805.1717, 2008.

\bibitem{Apo90}
Tom~M. Apostol.
\newblock {\em Modular Functions and Dirichlet Series in Number Theory}.
\newblock Springer, 2nd ed. edition, 1990.

\bibitem{Brocot1860}
A.~Brocot.
\newblock Calcul des rouages par approximation, nouvelle m{\'e}thode.
\newblock {\em Revue Chronom{\'e}trique}, 6:186--194, 1860.

\bibitem{Den38}
Arnaud Denjoy.
\newblock Sur une fonction r{\'e}elle de minkowski.
\newblock {\em Journal de Math{\'e}matiques Pures et Appliqu{\'e}es},
  17:105--151, 1938.

\bibitem{Dri99}
Dean~J. Driebe.
\newblock {\em Fully Chaotic Maps and Broken Time Symmetry}.
\newblock Kluwer Academic Publishers, 1999.

\bibitem{Dunford1957}
Nelson Dunford and Jacob~T. Schwartz.
\newblock {\em Linear Operators, Part I: General Theory}.
\newblock John Wiley \& Sons, 1957.
\newblock ISBN 0-471-60848-3.

\bibitem{Dushistova2007}
Anna~A. Dushistova and Nikolai~G. Moshchevitin.
\newblock On the derivative of the minkowski question mark function ?(x).
\newblock {\em ArXiv}, 0706.2219, 2007.

\bibitem{Gas92}
Pierre Gaspard.
\newblock r-adic one-dimensional maps and the euler summation formula.
\newblock {\em Journal of Physics A}, 25:L483--L485, 1992.

\bibitem{Glimm1981}
James Glimm and Arthur Jaffe.
\newblock {\em Quantum Physics, A Functional Integral Point of View}.
\newblock Springer-Verlag, 1981.
\newblock ISBN 0-387-90562-6.

\bibitem{Gut88}
Martin Gutzwiller and Benoit Mandelbrot.
\newblock Invariant multifractal measures in chaotic hamiltonian systems and
  related structures.
\newblock {\em Physical Review Letters}, 60:673--676, Feb 1988.

\bibitem{Gutzwiller1990}
Martin~C. Gutzwiller.
\newblock {\em Chaos in Classical and Quantum Mechanics}.
\newblock Springer-Verlag, 1990.
\newblock ISBN 0-387-97173-4.

\bibitem{Has92}
Hiroshi~H. Hasegawa and William C.Saphir.
\newblock Unitarity and irreversibility in chaotic systems.
\newblock {\em Physical Review A}, 46(12):7401--7423, 1992.

\bibitem{Iso03}
Stefano Isola.
\newblock On the spectrum of farey and gauss maps.
\newblock preprint, between 2000 and 2004.

\bibitem{Kac1959}
Mark Kac.
\newblock On the partition function of a one-dimensional gas.
\newblock {\em Physics of Fluids}, 2:8--12, 1959.

\bibitem{Kesse2008}
Marc Kesseb{\"o}hmer and Bernd~O. Stratmann.
\newblock Fractal analysis for sets of non-differentiability of minkowski's
  question mark function.
\newblock {\em Journal of Number Theory}, 128:2663--2686, 2008.

\bibitem{Khin35}
A.~Ya. Khinchin.
\newblock {\em Continued Fractions}.
\newblock Dover Publications, (reproduction of 1964 english translation of the
  original 1935 russian edition) edition, 1997.

\bibitem{Kittel1976}
Charles Kittel.
\newblock {\em Introduction to Solid State Physics}.
\newblock John Wiley \& Sons, 5th edition, 1976.
\newblock ISBN 0-471-49024-5.

\bibitem{Kle08}
Achim Klenke.
\newblock {\em Probability Theory}.
\newblock Springer, 2008.
\newblock ISBN 978-1-84800-047-6.

\bibitem{Man95}
B.B. Mandelbrot.
\newblock The minkowski measure and multifractal anomalies in invarient
  measures of parabolic dynamical systems.
\newblock In {\em Fractals and Chaos, The Mandelbrot Set and Beyond}, pages
  239--250. Springer, 2004.
\newblock reprint.

\bibitem{May91}
Dieter~H. Mayer.
\newblock Continued fractions and related transformations.
\newblock In C.~Series T.~Bedford, M.~Keane, editor, {\em Ergodic Theory,
  Symbolic Dynamics and Hyperbolic Spaces}, chapter~7, pages 175--222. Oxford
  University Press, 1991.

\bibitem{Mink04}
H.~Minkowski.
\newblock Zur geometrie der zahlen.
\newblock In {\em Verhandlungen des III. internationalen
  Mathematiker-Kongresses in Heidelberg}, pages 164--175. Berlin, 1904.

\bibitem{Paridis2001}
J.~Paridis, P.~Viader, and L.~Bibiloni.
\newblock The derivative of minkowski's ?(x) function.
\newblock {\em J. Math Anal. and Appl.}, 253:107--125, 2001.

\bibitem{Salem1943}
R.~Salem.
\newblock On some singular monotonic functions what are strictly increasing.
\newblock {\em Transactions of the American Mathematical Society}, 53:427--439,
  1943.

\bibitem{Stern1858}
M.~A. Stern.
\newblock {\"U}ber eine zehlentheoretische funktion.
\newblock {\em J. reine angew. Math.}, 55:193--220, 1858.

\bibitem{Ve-M04}
Linas Vepstas.
\newblock The minkowski question mark, gl(2,z) and the modular group.
\newblock http://www.linas.org/math/chap-minkowski.pdf, 2004.

\bibitem{Ve-T04}
Linas Vepstas.
\newblock Symmetries of period-doubling maps.
\newblock http://www.linas.org/math/chap-takagi.pdf, 2004.

\end{thebibliography}

\end{document}